\documentclass[11pt]{article}

\usepackage{amsmath, latexsym, amsfonts, amssymb, amsthm, amscd,graphicx, color
}
  \usepackage{graphics,epsf,psfrag}
  \newtheorem{theorem}{Theorem}[section]
  \newtheorem{lemma}[theorem]{Lemma}
  \newtheorem{proposition}[theorem]{Proposition}
  \newtheorem{corollary}[theorem]{Corollary}
  \newtheorem{remark}[theorem]{Remark}
  \newtheorem{definition}[theorem]{Definition}
    \newtheorem{condition}[theorem]{Condition}

\newcommand{\newatop}[2]{\genfrac{}{}{0pt}{}{#1}{#2}}

\def \P {{\mathbb P}}
\def \E {{\mathbb E}}

\newcommand{\cB}{{\cal B}}
\newcommand{\cG}{{\cal G}}
\newcommand{\cC}{{\cal C}}

\newcommand{\cM}{{\cal M}}
\newcommand{\cK}{{\cal K}}

\newcommand{\cS}{{\cal S}}

\newcommand{\cR}{{\cal R}}
\newcommand{\cW}{{\cal W}}

\newcommand{\cX}{{\cal X}}

\newcommand{\cI}{{\cal I}}

\newcommand{\cE}{{\cal E}}
\newcommand{\cF}{{\cal F}}

\def \gb {{\beta}}
\def \gep {{\varepsilon}}

\def \gl {{\lambda}}

\def \go {{\omega}}
\def \gd {{\delta}}

\newcommand{\be}[1]{\begin{equation}\label{#1}}
\newcommand{\ee}{\end{equation}}

\newcommand{\bl}[1]{\begin{lemma}\label{#1}}
\newcommand{\el}{\end{lemma}}

\newcommand{\br}[1]{\begin{remark}\label{#1}}
\newcommand{\er}{\end{remark}}

\newcommand{\bt}[1]{\begin{theorem}\label{#1}}
\newcommand{\et}{\end{theorem}}

\newcommand{\bd}[1]{\begin{definition}\label{#1}}
\newcommand{\ed}{\end{definition}}

\newcommand{\bcl}[1]{\begin{claim}\label{#1}}
\newcommand{\ecl}{\end{claim}}

\newcommand{\bp}[1]{\begin{proposition}\label{#1}}
\newcommand{\ep}{\end{proposition}}

\newcommand{\bc}[1]{\begin{corollary}\label{#1}}
\newcommand{\ec}{\end{corollary}}

\newcommand{\bi}{\begin{itemize}}
\newcommand{\ei}{\end{itemize}}

\newcommand{\ben}{\begin{enumerate}}
\newcommand{\een}{\end{enumerate}}

\def \qed {{\hfill$\square$\par\smallskip\par\noindent}}

\def \R {{\mathbb R}}

\def \N {{\mathbb N}}
\def \P {{\mathbb P}}
\def \E {{\mathbb E}}

\title{Metastability for general dynamics with rare transitions: escape 
time and critical configurations}

\author{
  E.N.M.\ Cirillo \footnote{ Dipartimento di Scienze di Base e Applicate per 
l'Ingegneria, 
 Sapienza Universit\`a di Roma, via A.\ Scarpa 16, I--00161, Roma, Italy. 
 \textit{emilio.cirillo@uniroma1.it}}\\
 F.R.\ Nardi  \footnote{ 
   Technische Universiteit Eindhoven, 
P.O. Box 513, 
5600 MB Eindhoven, The Netherlands.
 \textit{F.R.Nardi@tue.nl}}\\
 J.\ Sohier \footnote{ 
   Technische Universiteit Eindhoven, 
P.O. Box 513, 
5600 MB Eindhoven, The Netherlands.
\textit{j.sohier@tue.nl}   } }

\begin {document}
 

 \maketitle

   \begin{abstract}
Metastability is a physical phenomenon ubiquitous in first order 
phase transitions. 
A fruitful mathematical way to approach this phenomenon 
is the study of rare transitions Markov chains. 
For Metropolis chains associated with Statistical Mechanics systems, this phenomenon has been
described in an elegant
way in terms of the energy landscape associated 
to the Hamiltonian of the system. 
In this paper, we provide a similar description in the 
general rare transitions setup.
Beside their theoretical content, we believe that our results are a
useful tool to approach metastability for non--Metropolis systems 
such as Probabilistic Cellular Automata. 
\end{abstract}



%
\bigskip



%
%


\section{Introduction}
\label{s:intro}
In this paper we are interested in the phenomenon of metastability for systems evolving according to transformations 
     satisfying the thermodynamic law for small changes of the thermodynamical parameters. 
Metastability is a physical phenomenon ubiquitous in first order 
phase transitions. It is typically observed when a system is set up in a state which 
is not the most thermodynamically favored one and suddenly switches to the stable phase 
as a result of abrupt perturbations. 

Although metastable states have been deeply studied from the physical 
point of view, full rigorous mathematical theories based on a probabilistic
approach have been developed only in the last three decades. 
We refer to \cite{CNSo} for a complete recent bibliography. Let us just stress that the three main points of 
interest in the study of metastability
 are the description of: (i) the \textit{first hitting time} at which a Markov chain
  starting from the metastable state hits the stable one; (ii) the \textit{critical configurations} that the system
   has to pass to reach the stable states; (iii) the \textit{tube of typical trajectories} that the system 
    typically follows on its transition to the stable state. These notions are central quantities of interest 
   in many studies on metastability, which focus on proving convergence results
 in physically relevant limits, the most typical ones being the zero temperature limit and the infinite volume regime. 
In this paper, we focus on the finite volume and zero 
temperature limit setup.

The first mathematically rigorous results were obtained via the \emph{pathwise} approach, which has been first developed in 
the framework of special models and then fully understood in 
the context of the Metropolis dynamics \cite{CGOV,OS1,OV}.
In this framework, the properties of the first hitting time 
to the stable states are deduced via large deviation estimates 
on properly chosen tubes of trajectories. 
A different point of view, the \emph{potential theoretical} 
approach, has been proposed in \cite{BEGK04} and is based 
on capacity--like estimates. We mention that a more recent approach has also been developed in 
\cite{BL1,BL2}.

Here we adopt the pathwise point of view and generalize the 
theory to the general \emph{Freidlin--Wentzel Markov chains} 
or \emph{Markov chains with rare transitions} setup.
For Metropolis chains associated to Statistical Mechanics systems 
and reversible with respect 
to the associated Gibbs measure, the metastability phenomenon can be 
described in an elegant and physically satisfactory
way via the energy landscape associated 
with the Hamiltonian of the system \cite{OS1,OV}. 
In particular the time needed by the system 
to hit the stable state can be expressed in terms of the height of 
the most convenient path (that is the path with minimal energetic cost) that the system has to follow on its way along 
the energy landscape to the stable state. 
Moreover, the state of the system at the top of such a  path
is a \emph{gate} configuration in the sense that, in the low 
temperature regime,  
the system necessarily has to go through it before hitting the stable 
state. 
This description is very satisfactory
from the physical point of view since both the typical time 
that the system spends in the metastable state before switching to the 
stable one and the mechanism that produces this escape can be quantified purely through the energy landscape. Let us mention 
 that a simplified pathwise approach was proposed in \cite{MNOS}, where the authors disentangled the study of the first hitting 
  time from the study of the set of critical configurations and of the tube of the typical trajectories. 

In this paper we show that a similar physically remarkable 
description can be given in the 
general rare transitions (Freidlin--Wentzel) framework, when the invariant 
measure of the system is a priori not Gibbsian. 
In this setup the pathwise study of metastability has been approached 
with a different scheme in \cite{OS2}, where the 
physical relevant quantities describing the metastable state 
are computed via a renormalization procedure. 
Here we 
show that the strategy developed in \cite{MNOS} can be extended to this setup at the cost of a higher  
 complexity of techniques. A typical way of proceeding is to redefine 
 the height of a path in terms 
of the exponential weight of the transition probabilities and 
of a function, the \emph{virtual energy}, associated to the 
low temperature behavior of the invariant measure. 
In other words we reduce the pathwise study of metastability 
in the general rare transition case to the solution of a 
variational problem within the landscape induced by this notion 
of path height, using as a main tool the general cycle theory developed in \cite{Ca,CaCe}. We stress that, unlike the Metropolis case, 
 this procedure cannot be applied only from the detailed analysis of the set of optimal paths,
 and that a finer description of the cycle landscape is needed to 
  perform the analysis. 

Besides their theoretical content, the main motivation of our results 
has been to provide a
useful tool to approach metastability for a well known class of non--Metropolis systems, namely the Probabilistic Cellular Automata
\cite{CNP,CLRS}. 
Indeed, in this case, it is possible to write the virtual energy 
in a rather simple way and then solve the difficult variational 
problems in the induced landscape
\cite{CNS08p,CNS08j,CN03}. 

The technical difficulties that we had to overcome 
are rather evident: giving a satisfactory mathematical description of metastability in a context 
 where no Hamiltonian is available is a priori rather challenging. We overcame this difficulty using two key ideas. 

First idea. In the seminal papers on the pathwise approach to metastability 
\cite{OS1,OV} results were proved via detailed probability 
estimates on suitably chosen tube of trajectories. A simpler 
approach has been pointed out in \cite{MNOS}, where, 
still in the framework of the Metropolis dynamics, the author 
have shown that the main ingredient necessary to achieve the pathwise 
description of metastability is the classification 
of all the states of the systems in a sequence of 
decreasing (for the inclusion) subsets of the state space, whose elements have increasing stability,
in the sense that starting from any one of them the 
height that has to be bypassed to reach a lower energy level 
becomes increasingly higher. Moreover, the authors use in a crucial way 
a recurrence property stating that starting from 
any state, the process reaches one of these stability level sets 
within a time controlled exponentially by the stability level of 
the set itself. 
This is the point of view we also adopt in the present work. 

Second idea. 
One of the key tools in the pathwise study of metastability 
is the notion of cycle. In the context of general Markov 
chains, a cycle can be thought 
as a subset of the configuration states enjoying the following property: 
starting from anywhere 
within the cycle, with high probability the process visits all the 
states within 
the cycle before exiting the set itself. In the study of the metastable 
behavior of Metropolis chains a more physical definition of the notion of
cycle was used: a cycle is a set of configurations such that 
starting from any of them any other can be reached by a path within 
the set with maximal energy height smaller than the 
minimal one necessary for the process to exit the set. 
In this paper, following \cite{Ca}, we use the fact that by 
defining the height of a path in terms of the virtual energy and 
of the exponential cost of transition, the two different 
approaches to cycles can be proven to be equivalent. 

The paper is organized as follows. In Section~\ref{s:modello} we 
describe our setup and state the main results. 
Section~\ref{lecCa} is devoted to the discussion of the theory 
of cycles. 
In Section~\ref{s:pmr} we prove our main results. In Appendix \ref{s:gener}, we develop a condition under which the virtual energy 
 is explicitly computable, and in Appendix \ref{s:graph}, we make a quick recap about the virtual energy.

\section{Model and main results}\label{s:modello}
In this section we introduce a general setup 
and state our main results on the metastable behavior of such a system. Then we
describe in details this behavior in terms of the 
virtual energy, which in this setup is the analogous of the 
Hamiltonian for Metropolis chains. 

\subsection{The Freidlin--Wentzell setup}

 In this paper we will deal with a finite state space Markov chain with rare transitions. 
We consider
\begin{itemize}
\item[--]
an arbitrary 
finite \emph{state space} $\cX$.
\item[--]
A \emph{rate function} $\Delta: \cX\times\cX \mapsto \R^{+}\cup \{ \infty\}$. 
We assume that $\Delta$ is \emph{irreducible}
in the sense that 
for every $x,y \in \cX$, there exists a path 
$\go=(\go_{1},\ldots, \go_{n})\in\cX^n$ with $\go_{1} = x$, $\go_{n} = y$ and 
for every $ 1 \leq i \leq n-1, \Delta(x_{i},x_{i+1}) < \infty$, 
where $n$ is a positive integer.
\end{itemize}

\begin{definition}\label{setupFW}
A family of time homogeneous Markov chains 
$(X_{n})_{n \in \N}$ on $\cX$ with transition 
probabilities $p_{\beta}$ indexed by a positive parameter $\gb$
is said to 
\emph{"satisfy the Freidlin--Wentzell condition with respect to the rate 
function $\Delta$"} 
or 
\emph{"to have rare transitions with rate function $\Delta$"}
if and only if 
\begin{equation}
\label{FW}
\lim_{\gb \to \infty} \frac{-\log p_{\gb}(x,y)}{\gb} = \Delta(x,y)
\end{equation} 
for any $x,y\in\cX$.    
\end{definition}
     
The particular case where $\Delta(x,y)$ is infinite should be understood 
as the fact that, at low temperature, there is no transition possible 
between states $x$ and $y$. In many papers, a connectivity matrix is introduced, that is a matrix 
 whose non zero terms correspond to allowed jumps, see for instance
\cite{OV}[Condition R, Chapter~6].

We also note that condition \eqref{FW} is usually written explicitly; 
namely, for any $\gamma>0$, there exists $\beta_0>0$ such that 
\begin{equation}
\label{FW02}
e^{-\beta[\Delta(x,y)+\gamma]}
  \le p_\beta(x,y)
  \le e^{-\beta[\Delta(x,y)-\gamma]}
\end{equation}
for any $\beta>\beta_0$ and any $x,y\in\cX$. 
See for instance \cite{OV}[Condition~FW, Chapter~6] where 
the parameter $\gamma$ is assumed to be 
a function of $\beta$ vanishing for $\beta\to\infty$, so that in particular the Freidlin-Wentzell setup covers this case. 
    
\begin{remark}\label{partcas}
This framework covers in particular two relevant 
examples which have been under close scrutinity over the last decades.
\begin{enumerate}
\item \emph{Metropolis algorithm with Hamiltonian}
 $U:\cX \to \R$ (see, for instance, \cite{OV}[Condition~M, Chapter~6]
and \cite{MRRTT}). 
        It is the particular case where  
       \begin{equation}
       \label{metro}
       \Delta(x,y):=
       \left\{
       \begin{array}{ll}
       (U(y) - U(x))^{+} & \hspace{ 6 pt} \text{if} \hspace{6 pt} q(x,y)>0\\
       \infty & \hspace{6 pt} \text{otherwise}\\
       \end{array}
       \right.
       \;\;,
       \end{equation} 
       for any $(x,y) \in \cX\times\cX$ 
   where $q$ is an irreducible Markov matrix $\cX \times \cX \to [0,1]$ which does not depend on $\beta$. 
       We stress that the Metropolis algorithm itself is a general framework 
       which has as stationary measure the Gibbs measure of
       models issued from Statistical Mechanics 
       (see examples later). 
       \item \emph{Weak reversible dynamics with respect to 
        the potential} $U:\cX \to \R$ or 
        \emph{dynamics induced by the potential} 
        $U:\cX \to \R$.
  This is the case where the rate function $\Delta$ is 
  such that for any $(x,y) \in \cX\times\cX$
        \begin{equation}\label{weakrev}
         U(x) + \Delta(x,y) = U(y) + \Delta(y,x)
        \end{equation} 
 with the convention that $+ \infty + r = + \infty$ for any $r \in \mathbb{R}$. 
  
Even if the Metropolis dynamics is an example of a potential induced 
dynamics, these models form a broader class in which other important 
examples are 
Probabilistic Cellular Automata, see 
\cite{GJH,CNS08j,CN03} and the following Remark \ref{ossPCA}.  
\end{enumerate}
\end{remark}

From now on, we will always consider the general case 
of a family of homogeneous Markov chains 
satisfying
the condition in Definition~\ref{setupFW}.

\subsection{Virtual energy}
\label{virth}   
A fundamental notion for the physical 
approach of the problem of metastability in the setup of 
rare transitions chains is the notion 
of \emph{virtual energy}, whose definition is based on 
the following result. 
     
\begin{proposition}[\cite{Ca}[Proposition 4.1] 
\label{virtH}
Consider
a family of Markov chains satisfying the Freidlin--Wentzell 
condition in Definition~\ref{setupFW}.
For $\gb$ large enough, each Markov chain is 
irreducible and its invariant probability distribution 
$\mu_{\gb}$ is such that for 
          any $x \in \cX$, the limit 
\begin{displaymath}
\lim_{\gb \to \infty} -\frac{1}{\gb} \log\mu_{\gb}(x) 
\end{displaymath} 
exists and is a positive finite real number.
      \end{proposition}

\begin{definition}
\label{d:virth}
In view of Proposition~\ref{virtH}, 
the limiting function 
\begin{equation}
\label{ve000}
H(x):=\lim_{\gb \to \infty} -\frac{1}{\gb} \log\mu_{\gb}(x),
\end{equation}
for $x\in\cX$, 
is called \emph{virtual energy}.
\end{definition}

       The proof of Proposition \ref{virtH} relies on some deep combinatorial results which are tailored to the Freidlin--Wentzell 
        context. In general, the virtual energy has an exact expression 
in function of the transition rates $\Delta$ (see, for instance,
\cite{Ca}[Proposition~4.1], or the Appendix \ref{s:graph} at the end of the present work).
Unfortunately, in the most general setup, this expression involving a certain family of graphs 
is intractable for all practical purposes when one is interested 
to study particular models. 

 \begin{remark}\label{ossPCA}
  In the special case of 
Probabilistic Cellular Automata, \cite{CN03,CNS08j}, the authors deal 
 with  models involving a potential 
$G_{\gb}(x)$ depending on $\gb$ and satisfying the balance condition
\begin{displaymath}
p_{\gb}(x,y) e^{-G_{\gb}(x)} = p_{\gb}(y,x) e^{-G_{\gb}(y)} 
\end{displaymath} 
for every positive $\gb$. 
To bypass the technical difficulties inherent to these models, 
which stem for a large part
from the intricate dependence 
on $\gb$ of $p_{\gb}(\cdot)$ and $G_{\gb}(\cdot)$, 
the authors computed directly the expressions of the  
rate function $\Delta(\cdot)$ in \eqref{FW}
and of the virtual energy \eqref{ve000}. In this way, they obtained a weak reversible dynamics (see \eqref{weakrev}). 
It thus became easier to solve 
the metastable behavior 
for these models, using solely the limit expressions obtained. 
We refer to Appendix \ref{s:gener} for a more general 
context in which these techniques still apply and we mention that 
our hope is that this generalization should cover some 
other relevant cases in which only the transitions rates are explicitly 
computable. 
  
 \end{remark}

Finally, we stress that in the particular 
cases of Remark~\ref{partcas}, the virtual energy, up to an additive constant, 
is precisely the potential which induces the dynamics. 
          
\begin{proposition}[\cite{Ca}Proposition~4.1]
\label{t:ve020}
In the particular case of the dynamics induced by the 
potential $U:\cX\to\R$
(see Remark~\ref{partcas}) one can show the equality
\begin{displaymath}
H(x) = U(x) - \min_{\cX} U
\end{displaymath} 
for any $x\in\cX$.
\end{proposition}
     
\subsection{General definitions}
\label{s:defpre}

In the present and in the following sections, we 
introduce some standard notions, which are natural generalizations of the 
analogous quantities in the reversible setup, see \cite{MNOS} or \cite{OV}. 

A real valued function $f: \R^{+} \mapsto \R^{+}$ is \emph{super exponentially 
small} (SES for short) if and only if
      \begin{displaymath}
       \lim_{\gb \to \infty} \frac{1}{\gb}\log f(\gb) = - \infty.
      \end{displaymath} 

For $x \in \cX$, we let $X^x_t$ 
be the chain started at $x$.
For a nonempty set $A\subset\cX$ and $x \in \cX$, we 
introduce
the \emph{first hitting time} $\tau_{A}^{x}$ to the set $A$ 
which is 
the random variable 
       \begin{displaymath}
        \tau_{A}^{x}: = \inf \{ k \geq 0, X_{k}^{x} \in A \}.
       \end{displaymath} 

A \emph{path} is a sequence $\go= (\go_{1},\ldots, \go_{n})$
       such that $\Delta(\go_{i},\go_{i+1}) < \infty$ for $i=1,\ldots, n-1$.
For a path $\go= (\go_{1},\ldots, \go_{n})$, we define $|\go| = n$ its 
\emph{length}. 
For $x,y\in\cX$ 
a path $\go: x \to y$ \emph{joining $x$ to $y$} is a path 
$\go= (\go_{1},\ldots, \go_{n})$
such that $\go_{1} = x$ and $\go_{n} = y$. 
For any $x,y\in\cX$
we write $\Omega_{x,y}$ for the set of paths joining $x$ to $y$.
For $A,B\subset\cX$ nonempty sets, we write 
$\Omega_{A,B}$ for the set of paths joining a point in $A$ to 
a point in $B$.

A set $A \subset \cX$ with $|A| > 1$ is \emph{connected} if and only 
if for all $x,y\in A$, there exists a path 
          $\go\in\Omega_{x,y}$ such that for any 
$i \leq |\go|, \go_{i} \in A$. 
By convention, we say that every singleton is connected. 

For a nonempty set $A$, we define its \emph{external boundary} 
$\partial A:=\{y\in\cX\setminus A, \text{there exists } x \in A
                 \textrm{ such that } \Delta(x,y) < \infty \}$
and we write 
\begin{equation}
\label{fondo}
H(A) = \min_{A}H.
\end{equation}

The \emph{bottom} $\cF(A)$ of $A$ is the set of global minima
of $H$ on $A$, that is
\begin{displaymath}
\cF(A): = \{x \in A, H(x) = H(A) \}.
\end{displaymath} 
The set $\cX^\textrm{s}:=\cF(\cX)$ is called the set of 
\emph{stable points} or the set of \emph{ground states} 
of the virtual energy.

\subsection{Communication height}
\label{s:ch}
A key notion in studying metastability is the one of the cost 
that the chain has to pay to follow a path. 
In the case of Metropolis dynamics this quantity is the highest 
energy level reached along a path. Such a notion has to be modified 
when general rare transitions dynamics are considered \cite{Tr,CN03}. 
We thus define the
\emph{height} or \emph{elevation} $\Phi(\omega)$ of a path 
$\omega = (\omega_{1}, \ldots, \omega_{n})$ by setting 
\begin{equation}
\label{ch000}
\Phi(\omega):= \max_{i=1,\ldots,|\go|-1} [H(\omega_{i}) 
                              + \Delta(\omega_{i},\omega_{i+1})].
\end{equation} 
 
The \emph{communication height} $\Phi(x,y)$ between two states 
$x,y \in \cX$ is the quantity
\begin{equation}
\label{ch010}
\Phi(x,y):= \min_{\omega \in \Omega_{x,y}} \Phi(\omega).
\end{equation} 
 
Given two nonempty sets $A,B\subset\cX$, we define 
\begin{equation}
\label{optpa}
\Phi(A,B):= \min_{x \in A, y \in B} \Phi(x,y)
\end{equation}          
the \emph{communication height} between $A$ and $B$. 

For $A,B$ nonempty subsets of $\cX$, we define 
$\Omega_{A,B}^\textrm{opt}$ as the set of \emph{optimal paths joining 
$A$ to $B$}, that is the 
set of paths joining a point in $A$ to a point in $B$ and realizing the 
min--max $\Phi(A,B)$ defined in \eqref{optpa}. 

For rare transitions dynamics induced by a potential 
(see Remark~\ref{partcas}) it is easy to see 
that the communication height between two states 
is symmetric. 
A non--trivial result due to A.\ Trouv\'e \cite{Tr}
states that this is the case even in the general setup adopted in 
this paper. 

\begin{proposition}[\cite{Ca} Proposition 4.14]
\label{sym}
The communication height between states is symmetric, that is, 
$\Phi(x,y) = \Phi(y,x)$ for any $x,y\in\cX$.
\end{proposition}

\begin{corollary}[\cite{Ca} Proposition~4.17 ]
\label{t:notsym}
For any $x,y\in\cX$, the virtual energy satisfies
\begin{displaymath}
H(y)\le H(x)+\Delta(x,y).
\end{displaymath}
\end{corollary}

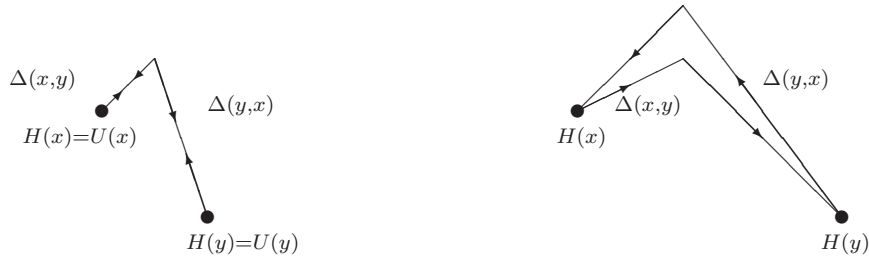
\begin{figure}
\begin{picture}(200,100)(5,15)
\put(100,60){\circle*{5}}
\put(140,20){\circle*{5}}
\put(100,60){\line(1,1){20}}
\put(100,60){\vector(1,1){8}}
\put(120,80){\vector(-1,-1){8}}
\put(140,20){\line(-1,3){20}}
\put(140,20){\vector(-1,3){8}}
\put(120,80){\vector(1,-3){8}}
\put(70,48){$\scriptstyle H(x)=U(x)$}
\put(132,8){$\scriptstyle H(y)=U(y)$}
\put(65,70){$\scriptstyle \Delta(x,y)$}
\put(140,60){$\scriptstyle \Delta(y,x)$}
\put(280,60){\circle*{5}}
\put(380,20){\circle*{5}}
\put(280,60){\line(2,1){40}}
\put(280,60){\vector(2,1){20}}
\put(380,20){\line(-1,1){60}}
\put(320,80){\vector(1,-1){30}}
\put(272,48){$\scriptstyle H(x)$}
\put(372,8){$\scriptstyle H(y)$}
\put(294,60){$\scriptstyle \Delta(x,y)$}
\put(280,60){\line(1,1){40}}
\put(320,100){\vector(-1,-1){20}}
\put(320,100){\line(3,-4){60}}
\put(380,20){\vector(-3,4){40}}
\put(350,70){$\scriptstyle \Delta(y,x)$}
\end{picture}
\vskip 0.5 cm
\caption{Illustration of the result in Corollary~\ref{t:notsym}. The 
picture on the left refers to the weak reversible case, 
whereas the picture on the right refers to the general dynamics with rare 
transitions.}
\label{f:coro}
\end{figure}

This corollary is quite interesting and its meaning is 
illustrated in Figure~\ref{f:coro}. Indeed, in the 
case of a dynamics induced by a potential, the jump between 
two states can be thought of as in the left part of the figure: 
the chain can jump in both directions and the height reached 
in both cases is the same. 
This is not true anymore in general under the sole 
assumptions of Definition~\ref{setupFW} (see the illustration on
the right in the same figure). Provided the chain can perform 
the jump from $x$ to $y$, that is $\Delta(x,y)<\infty$, it is not 
ensured that the reverse jump is allowed. Moreover, even in such 
a case, the heights which are attained during the two jumps 
in general are different. 
Nevertheless, the important Corollary~\ref{t:notsym}
states that the virtual energies of the two states $x$ and $y$ 
are both smaller than the heights attained by performing 
any of the two jumps.

\subsection{Metastable states}
\label{defmet}
The main purpose of this article is to define the notion of metastable 
states for a general rare transition dynamics and to prove 
estimates on the hitting time to the set of stable states for the 
dynamics started at a metastable state. 

To perform this, we need to introduce the notion of 
stability level of a state $x \in \cX$. First define 
\begin{equation}
\label{sl000}
\cI_{x}: = \{ y \in \cX, H(y) < H(x) \}
\end{equation} 
which may be empty in general. 
Then we define the \emph{stability level} of any state $x\in\cX$ by
\begin{equation}
\label{sl001}
V_{x}: = \Phi(x,\cI_{x}) - H(x)
\end{equation} 
and we set $V_{x} = \infty$ in the case where $\cI_{x}$ is empty. 
We also let 
\begin{equation}
\label{sl002}
V^\textrm{m}: = \max_{x\in\cX\setminus\cX^\textrm{s}}V_{x}
\end{equation} 
be the \emph{maximal} stability level.
       
Metastable states should be thought of as the set of states 
where the dynamics is typically going to spend a lot of time before 
reaching in a drastic way the set of stable states $\cX^\textrm{s}$.
Following \cite{MNOS} we define 
the set of \emph{metastable states} $\cX^\textrm{m}$
as 
\begin{equation}
\label{meta000}
\cX^\textrm{m}
:= 
\{x\in\cX,\, V_{x} = V^\textrm{m}\}
\end{equation} 
and in the sequel, see Section~\ref{s:meta}, we will state some 
results explaining why $\cX^\textrm{m}$ meets the 
requirements that one would heuristically expect from the set of 
metastable states. For example, we prove that the maximal stability 
level $V^\textrm{m}$ is precisely the quantity controlling
the typical time that the system needs to escape from the metastable 
state. 

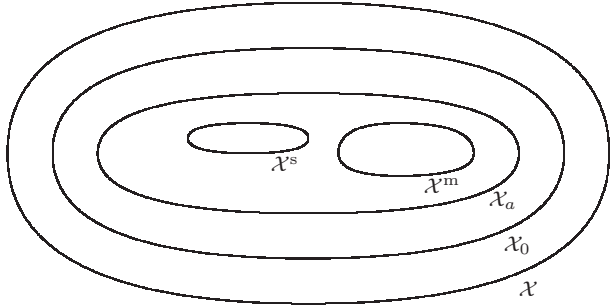
\begin{figure}[t]
 \centering
 \begin{picture}(400,150)(-60,-15)
 \setlength{\unitlength}{0.02cm}
 \thinlines
 \qbezier(0,100)(0,200)(200,200)
 \qbezier(0,100)(0,0)(200,0)
 \qbezier(200,200)(400,200)(400,100)
 \qbezier(200,0)(400,0)(400,100)
 \qbezier(30,100)(30,170)(200,170)
 \qbezier(30,100)(30,30)(200,30)
 \qbezier(200,170)(370,170)(370,100)
 \qbezier(200,30)(370,30)(370,100)
 \qbezier(60,100)(60,140)(200,140)
 \qbezier(60,100)(60,60)(200,60)
 \qbezier(200,140)(340,140)(340,100)
 \qbezier(200,60)(340,60)(340,100)
 \qbezier(220,100)(220,120)(260,120)
 \qbezier(220,100)(220,85)(260,85)
 \qbezier(260,120)(310,120)(310,100)
 \qbezier(260,85)(310,85)(310,100)
 \qbezier(120,110)(120,120)(160,120)
 \qbezier(120,110)(120,100)(160,100)
 \qbezier(160,120)(200,120)(200,110)
 \qbezier(160,100)(200,100)(200,110)
 \put(340,5){${\scriptstyle \cX}$}
 \put(330,35){${\scriptstyle \cX_0}$}
 \put(320,65){${\scriptstyle \cX_a}$}
 \put(277,74){${\scriptstyle \cX^{\textrm{m}}}$}
 \put(175,88){${\scriptstyle \cX^{\textrm{s}}}$}
 \end{picture}
\caption{Illustration of the structure of the 
          sets $\cX_a$'s (see definition \eqref{mestru})
         with $0<a<V^\textrm{m}$.}
\label{f:mestru}
\end{figure}

More generally, for any $a>0$,
we define the \emph{metastable set of level $a>0$} 
as follows
\begin{equation}
\label{mestru}
\cX_{a}:= \{ x \in \cX, V_{x} > a\}.
\end{equation} 
The structure of the sets $\cX_a$'s is depicted in 
Figure~\ref{f:mestru}.
It is immediate to realize that $\cX_a\subset\cX_{a'}$ for $a\ge a'$. 
Moreover, it is worth noting that 
$\cX_{V^\textrm{m}}=\cX^\textrm{s}$.

\subsection{Cycles, saddles, and gates}
\label{s:sg}   
%

 We stress that one of our main results (see Theorem \ref{crossgate} below) describes a family of sets which have to be crossed with large probability 
  in the low temperature limit. 
   
To introduce these sets, we define as in \cite{MNOS} the notion of \textit{saddle points} and of \textit{gates}. 
  We stress that, unlike the Metropolis dynamics, these notions cannot be defined at the level of paths only. Let us discuss this point 
   a bit since this is a major difference between the setups.

  The following definition was introduced in \cite{MNOS}, and we recall it for expository purposes only. We stress that we cannot 
   adapt it straightforwardly to our setup, as is discussed below.
   
   It would be natural to generalize the definition (see \cite{MNOS}) of the set of \emph{minimal saddles} 
between two states $x,y\in\cX$  in 
the context of Metropolis dynamics as
 \begin{displaymath}
  \begin{aligned}
\big \{z \in \cX, &\text{ there exists } 
                    \go\in\Omega_{x,y}^\textrm{opt}
                    \text{ and } i \leq |\go|\\
         & \text{ such that } \go_{i} = z 
       \text{ and } H(\go_{i-1}) + 
      \Delta(\go_{i-1},\go_{i}) = \Phi(x,y)  \big \}.
  \end{aligned}
 \end{displaymath} 
 
 In the  Freidlin Wentzell setup, this precise definition does not make sense at the level of typical 
     behavior of trajectories. For example, there might
     be an optimal path $\omega$ joining $x$ to $y$ and a 
      a minimal gate $W$ such that $\omega_{i} \in W$ (and hence  $H(\omega_{i-1}) + \Delta(\omega_{i-1},\omega_{i}) = \Phi(x,y)$)
      and such that nevertheless the point $\omega_{i}$  does not play any particular role for the dynamics. Indeed,
      there might be a path with cost strictly lower than $\Phi(x,y)$ joining 
        $\omega_{i-1}$ to $\omega_{i}$ which will be favoured by the dynamics in the low temperature limit.
        
         This phenomenon is very peculiar to the Metropolis setup; indeed, an energy level has to correspond to a point in this setup, whereas 
         in the Freidlin Wentzell setup, this correspondence is not valid anymore. 
         
         Nevertheless, we stress that we can  generalize 
         the notion of gates and of minimal gates in our setup, at the cost of higher complexity of definitions. To perform this, we need 
         to introduce the key 
          notions of \textit{cycle} 
and of \textit{principal boundary of a set}. 
The notion of cycle will be discussed 
in details in Section~\ref{lecCa}.
          
          \begin{definition}[\cite{Ca} Definition~4.2 ] 
\label{defcy}
A nonempty set $C\subset\cX$ is a \emph{cycle} if it is either a 
singleton or for any $x,y\in C$, such that $x\neq y$, 
\begin{displaymath}
\lim_{\beta\to\infty}-\frac{1}{\gb}
       \log\P_\gb[X_{\tau^x_{(\cX\setminus C)\cup\{y\}}}\neq y]>0.
\end{displaymath}
\end{definition}
 
In words, 
a nonempty set $C\subset\cX$ is a \emph{cycle} if it is either a 
singleton or for any $x,y\in C$ such that $x\neq y$, the 
probability starting from $x$ to leave $C$ without visiting $y$ is 
exponentially small.
We will denote by $\cC(\cX)$ the set of cycles. 
The set  $\cC(\cX)$ has a tree structure, that is:

\begin{proposition}[\cite{Ca}[Proposition~4.4]
\label{tree}
For any pair of cycles $C,C'$ such that $C \cap C' \neq \emptyset$, 
either $C \subset C'$ or $C' \subset C$. 
\end{proposition}

Next we introduce the important notion of \textit{principal boundary} 
of an arbitrary subset of the state space $\cX$.  
 
\begin{proposition}[\cite{Ca} Proposition~4.2]
\label{exgen}
 For any $D \subset \cX$ and any $x \in D$, the following limits exist
and are finite:
\begin{equation}
\label{rafuE}
\lim_{\gb \to \infty} 
\frac{1}{\gb} \log\E_{\gb}[ \tau^x_{\cX\setminus D}] 
=:\Gamma_{D}(x)
\end{equation} 
and,
for any $y\in\cX\setminus D$,
\begin{equation}
\label{rafu}
\lim_{\gb \to \infty} 
-\frac{1}{\gb} \log\P_{\gb}[X_{\tau^x_{\cX\setminus D}} = y] 
=:\Delta_{D}(x,y).
\end{equation} 
\end{proposition}

We stress that the limits appearing in the right hand side of 
\eqref{rafu} and \eqref{rafuE} have explicit expressions which, as in Definition \ref{d:virth} for the virtual energy, 
seem to be intractable for practical purposes at least in the field of statistical mechanics. 

The meaning of the two functions introduced in the Proposition \ref{exgen}
is rather transparent: \eqref{rafuE} provides an exponential 
control on the typical time 
needed to escape from a general domain $D$ starting from 
a state $x$ in its interior and $\Gamma_D(x)$ is the mass of such 
an exponential. On the other hand, \eqref{rafu} provides 
an exponential bound to the probability to escape from $D$, starting 
at $x$, through the site $y\in\cX\setminus D$. Hence, we can think 
to $\Delta_D(x,y)$ as a measure of the cost that has to be paid 
to exit from $D$ through $y$. 

Now, we remark that,
due to the fact that the state space $\cX$ is finite, 
for any domain $D\subset\cX$ and for any $x\in D$ 
there exists at least a point 
$y\in\cX\setminus D$ such that $\Delta_D(x,y)=0$. Thus, we 
can introduce the concept of \emph{principal boundary of a set $D\subset\cX$}
\begin{equation}\label{defBo}
\cB(D):=\{ y \in\cX\setminus D,\, \Delta_{D}(x,y) = 0 
           \text{ for some } x \in D \}. 
\end{equation}

 
  We are finally ready to describe in a rigorous way the notion of \textit{gates} which will be used to state one of our main results,  
   Theorem \ref{crossgate}. 
 
   \begin{definition}\label{defgateFW}
    Let $x,y \in \cX$. Let $C_{x,y}$ be the minimal cycle containing both $x$ and $y$  and let $\cM_{x,y} = \{C_{i}, i \leq n_{0} \}$
    be its decomposition into maximal strict subcycles. Both these notions are well defined by Proposition \ref{tree}. We define the set 
      of \textit{saddles} 
between $x$ and $y$ (which is denoted by $\cS(x,y)$) by
       \begin{displaymath}
        \cS(x,y) = \bigcup_{C \in \cM_{x,y}} \cB(C).
       \end{displaymath}

   \end{definition}

   We stress that the set $\cS(x,y)$ is related in a very intricate way to the energy landscape of the dynamics. 
   
   From now on, we can proceed by analogy with the definitions of the Metropolis case (see \cite{MNOS}). 
Given $x,y \in \cX$, we say that 
$W \subset\cX$ is a \emph{gate for the couple $(x,y)$} 
if 
$W \subset \cS(x,y)$ and every path in $\Omega_{x,y}^\textrm{opt}$ 
intersects $W$, that is 
\begin{displaymath}
\go \in\Omega_{x,y}^\textrm{opt} \Longrightarrow 
\go \cap W \neq \emptyset.
\end{displaymath} 
We also introduce $\cW(x,y)$ as being the collection of all the gates 
for the couple $(x,y)$. 

A gate $W\in\cW(x,y)$ for the $(x,y) \in\cX$ is \emph{minimal} 
if it is a minimal (for the inclusion relation) element of 
$\cW(x,y)$. Otherwise stated, for 
any $W'\subsetneq W$, 
there exists $\go'\in\Omega_{x,y}^\textrm{opt}$ such that 
$\go'\cap W'=\emptyset$. 
In the metastability literature, the following  set is also standard 
  \begin{displaymath}
  \cG(x,y) := \bigcup_{W \in \cW(x,y), W \text{ is minimal.} } W;
  \end{displaymath}      
namely, $\cG(x,y)$ is  the set of saddles between $x$ and $y$ belonging 
to a minimal gate in $\cW(x,y)$.
   
%

 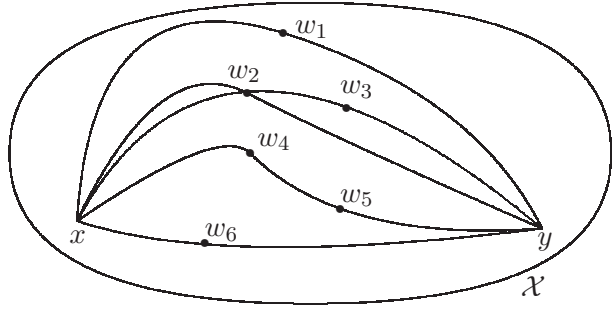
\begin{figure}[t]
  \centering
  \begin{picture}(400,150)(-60,-15)
  \setlength{\unitlength}{0.02cm}
  \thinlines
  \qbezier(0,100)(0,200)(200,200)
  \qbezier(0,100)(0,0)(200,0)
  \qbezier(200,200)(400,200)(400,100)
  \qbezier(200,0)(400,0)(400,100)
  \qbezier(45,55)(55,220)(182,180)
  \qbezier(182,180)(309,140)(354,50)
  \put(182,180){\circle*{5}}
  \put(190,180){$w_1$}
 
  \qbezier(45,55)(100,170)(158,140)
  \qbezier(158,140)(216,110)(354,50)
  \put(158,140){\circle*{5}}
  \put(145,150){$w_2$}
 
  \qbezier(45,55)(140,230)(354,50)
  \put(224,130){\circle*{5}}
  \put(220,140){$w_3$}
 
  \qbezier(45,55)(140,120)(160,100)
  \qbezier(160,100)(220,40)(354,50)
  \put(160,100){\circle*{5}}
  \put(165,105){$w_4$}
  \put(220,63){\circle*{5}}
  \put(220,68){$w_5$}
  
  \qbezier(45,55)(140,23)(354,50)
  \put(130,40){\circle*{5}}
  \put(130,45){$w_6$}
  \put(340,5){$\cX$}
  \put(40,40){$x$}
  \put(351,39){$y$}
  \end{picture}
 \caption{Illustration of the notion of gate between two 
          configurations $x$ and $y$. The case depicted here 
          is the following:
          $\cS(x,y)=\{w_1,\ldots,w_6\}$. The optimal paths in 
          $\Omega^\textrm{opt}_{x,y}$ are represented by the 
          five black lines.
          The minimal gates are $\{w_1,w_2,w_4,w_6\}$ and 
          $\{w_1,w_2,w_5,w_6\}$.
          Any other subset of $\cS(x,y)$ obtained by adding 
          some of the missing saddles to one of the two 
          minimal gates is a gate. 
 }
 \label{f:gate}
 \end{figure}
 
%

%
\subsection{Main results}
\label{s:meta}
In this section we collect our results about the behavior of 
the system started at a metastable state. These results justify a posteriori why the abstract notion of metastable set 
$\cX^\textrm{m}$ fits with the heuristic idea of 
metastable behavior.

The first two results state that the escape time, that is the 
typical time needed by the dynamics started at a metastable state 
to reach the set of stable states, is exponentially large in 
the parameter $\beta$. Moreover, they ensure that 
the mass of such an exponential is given by the 
maximal stability level; the first result is a convergence 
in probability, whereas the second ensures convergence in mean.

\begin{theorem}
\label{t:cp}
For any $x \in \cX^{\normalfont\textrm{m}}$, 
for any $\varepsilon>0$ there exists $\beta_0<\infty$ and 
$K>0$ such that 
\begin{equation}
\label{cp01}
\P_\gb[
\tau^{x}_{\cX^{\normalfont\textrm{s}}}
<
e^{\beta(V^{\normalfont\textrm{m}}-\varepsilon)}
]
<e^{-\beta K}
\end{equation}
and 
\begin{equation}
\label{cp02}
{\normalfont\textit{ the function }}
\beta\mapsto
\P_\gb[
\tau^{x}_{\cX^{\normalfont\textrm{s}}}
>
e^{\beta(V^{\normalfont\textrm{m}}+\varepsilon)}
]
{\normalfont \;\;\textit{ is SES}}.
\end{equation}
\end{theorem}

\begin{theorem} 
\label{tune}
For any $x \in \cX^{\normalfont\textrm{m}}$, 
the following convergence holds 
\begin{equation}
\label{e:tune}
\lim_{\gb\to\infty}\frac{1}{\gb} 
  \log \E_\beta[\tau^{x}_{\cX^{\normalfont\textrm{s}}}]
=
V^{\normalfont\textrm{m}}.
\end{equation} 
\end{theorem}
  
\begin{theorem}
\label{expbeh} 
Assume the existence of a recurrent state $x_{0}$ for the dynamics,
namely, assume that there exists $x_{0} \in \cX$ 
   such that
    \begin{itemize}
     \item[--] \textbf{ late escape from the state $x_{0}$:} 
      \begin{equation}\label{le}
       T_{\gb} := \inf \left \{ n \geq 0, 
                \P_\gb 
       \left[\tau_{\cX^{\normalfont\textrm{s}}}^{x_{0}}\le n\right] 
                     \geq 1 - e^{-1} \right \} 
              \stackrel{\gb\to\infty}{\longrightarrow} \infty;
      \end{equation}
       
   \item[--] \textbf{ fast recurrence to $x_{0}$:}
   there exist two functions 
   $\delta_\beta,T'_\beta:[0,+\infty]\to\R$
   such that 
   \begin{equation}
   \label{fr0}
   \lim_{\gb \to \infty} \frac{T_{\gb}'}{T_{\gb}} = 0 ,\;\;\;
   \lim_{\gb \to \infty} \gd_{\gb} = 0, 
   \end{equation}
   and
   \begin{equation}
   \label{fr}
   \P_\gb \Big[\tau_{\{x_{0},\cX^{\normalfont\textrm{s}}\}}^{x} 
          > T_{\gb}' \Big] \leq \gd_{\gb}
   \end{equation}
   for any $x\in\cX$ and $\beta$ large enough. 
   \end{itemize}
  
\par\noindent
Then, the following holds
      \begin{enumerate}
      \item\label{i:expbeh1}
      the random variable 
       $\tau_{\cX^{\normalfont\textrm{s}}}^{x_{0}}/T_{\gb}$ converges in law to 
      an exponential variable with mean one;
      \item\label{i:expbeh2}
      the mean hitting time and $T_{\gb}$ are asymptotically equivalent, 
      that is
       \begin{equation}
       \label{eq:expbeh}
        \lim_{\gb \to \infty} \frac{1}{T_{\gb}}\,
                \E_\gb[\tau_{\cX^{\normalfont\textrm{s}}}^{x_{0}}] = 1;
       \end{equation} 
      \item\label{i:expbeh3}
      the random variable
      $\tau_{\cX^{\normalfont\textrm{s}}}^{x_{0}}
       /\E_\gb[\tau_{\cX^{\normalfont\textrm{s}}}^{x_{0}}]$ 
      converges in law to an exponential variable with mean one. 
      \end{enumerate}
  \end{theorem}
  
We stress that such exponential behaviors are not new in the literature; 
for the Metropolis case,
we refer of course to \cite[Theorem~4.15]{MNOS}, and we refer 
to \cite{AB1,AB2} for the generic reversible case. 
In an irreversible setup, results appeared only much more 
recently; let us mention \cite{BLM} and \cite{OI}. 
In the case where the cardinality of the 
state space $\cX$ diverges, 
more precise results than the one described in 
Theorem~\ref{expbeh} were obtained in \cite{FMNS} and \cite{FMNSS}. 

%
           
Our result is different from the ones we mention here, 
since we are able 
to give the explicit value of the expected value of the escape
time in function of the transition rates of the family of Markov chains. 

The above results are related to the properties of the escape 
time, the following one gives in particular some information about the trajectory that 
the dynamics started at a metastable state 
 follows with high probability on its way towards the stable state.

\begin{theorem}
\label{crossgate}
For any pair $x,y\in\cX$ we consider the set of gates $\cW(x,y)$
introduced in Section \ref{s:sg} and the corresponding set of minimal gates.
 For any minimal gate $W\in\cW(x,y)$, 
there exists $c>0$ such that 
\begin{displaymath}
\P_\gb[\tau_{W}^{x} > \tau_{y}^{x}] \leq e^{-\gb c}
\end{displaymath} 
for $\gb$ sufficiently large.
\end{theorem}

The typical example of application of this result is to consider $x\in\cX^\textrm{m}$, $y\in\cX^\textrm{s}$, and 
$W\in\cW(x,y)$; 
Theorem~\ref{crossgate} ensures that, with high probability, 
on its escape from the metastable state $x$, the dynamics has to visit 
the gate $W$ before hitting the stable state $y$. This is 
a strong information about the way in which the dynamics 
 performs its escape from a metastable state.

   We stress that our main tool to prove Theorem \ref{crossgate}  is the description in great details of the set of typical trajectories of the dynamics
    of the transitions from $x$ to $y$, which is the \textit{tube of typical trajectories} $\cK_{x,y}$ (see \cite[Chapter 6]{OV}, and in particular 
     Part 6.7, Theorems 6.31 and 6.33 where an analogous description has been performed in the particular case of the Metropolis dynamics). Recall the 
     notations $C_{x,y}$ (the minimal cycle containing $x$ and $y$) and $\cM_{x,y}$ (the decomposition into maximal strict subcycles of 
      $C_{x,y}$) of Definition \ref{defgateFW}. The set  $\cK_{x,y}$
     is a subset of $\Omega_{x,y}^{\rm{opt}}$ which can be described 
    as follows: 
    \begin{enumerate}
     \item as soon as the dynamics enters an element $C \in \cM_{x,y}$, it exits $C$ through its principal boundary  $\cB(C)$. 
     This implies in particular the fact that the dynamics stays within the cycle $C_{x,y}$ during its transition 
       from $x$ to $y$, as we will show later (see in particular Remark \ref{selinC}); 
      \item as soon as the dynamics enters the unique element $C(y)$ of $\cM_{x,y}$ containing $y$, it hits 
       $y$ before leaving $C(y)$ for the first time. 
     \end{enumerate}
     
     We then state the following proposition about the tube $\cK_{x,y}$:   
    \begin{proposition}\label{traje}
       For any $x,y \in \cX$, as $\gb \to \infty$, the set $\cK_{x,y}$ has probability exponentially close to $1$, that is, for 
        any $\gep > 0$, there exists $\gb_{0}$ such that for any $\gb \geq \gb_{0}$:
       \begin{displaymath}
        \P_{\gb} [\cK_{x,y}] \geq 1 - e^{-\gb \gep}. 
       \end{displaymath} 
    \end{proposition}
    
     We stress that in concrete models, such a detailed description of the exit tube relies on an exhaustive analysis of the energy 
      landscape which is unlikely to be performed in general. Nevertheless, for the particular case of PCA's, this 
  analysis can be greatly simplified.

  \begin{remark} \label{decsottocrit} 
     For reversible PCA's, the analysis of the 
    phenomenon of metastability was performed in \cite{CNS08j} by studying the transition between the metastable
    state (the $-$ phase) towards the stable 
     state (the $+$ phase in this specific model) using a particular case of Proposition \ref{traje}. 
     Indeed, the decomposition into maximal cycles $C_{(-),(+)}$ was reduced to two cycles only, and the one containing the $(-)$ 
      state 
     was refered to as the \textit{subcritical phase}.  One of the main tasks was then to identify the set of saddles, which 
     in this case was reduced 
     to the principal boundary of the subcritical phase. 
      
      Our approach shows in which way this technique should be extended in the more general case of several maximal cycles involved 
       in the maximal decomposition of the cycle $C_{x,y}$.  A practical way to perform this would be to use Definition 
       \ref{defgateFW} to identify recursively the set 
        of saddles.

  \end{remark}

\subsection{Further results on the typical behavior of trajectories}
\label{s:colla}
In this section we collect some results on the set of typical trajectories in the large $\gb$ limit. 
   
The first result of this section is a large deviation estimate 
on the hitting time to the metastable set $\cX_{a}$ at level $a>0$. 
The structure of the sets $\cX_a$'s is depicted in Figure~\ref{f:mestru}.
Given $a>0$, since states outside $\cX_a$ have stability level 
smaller that $a$, it is rather natural to expect that, starting 
from such a set, the system will typically need a time smaller than 
$\exp\{\beta a\}$ to reach $\cX_a$. This recurrence 
result is the content of the following lemma. 
      
\begin{proposition}
\label{Main}
For any $a>0$ and any $\gep>0$, the function 
\begin{displaymath}
\gb \mapsto  \sup_{x \in \cX} 
\P_\gb\big[ \tau_{\cX_{a}}^{x} > e^{\gb(a+ \gep) }\big] 
\end{displaymath} 
is SES. 
\end{proposition}

 \begin{remark} Proposition~\ref{Main} allows to disentangle the study of the first hitting time of the stable state from the results 
 on the tube of typical trajectories performed in great details both  in \cite{OS2} and in \cite{CaCe}. This remarkable fact 
  relies on Proposition \ref{grtyp}, which guarantees the existence
   of downhill cycle paths to exit from any given set. In the Metropolis setup, this has been performed in \cite{MNOS} (see Theorem 3.1 and 
    Lemma 2.28). 
   
 \end{remark}



The following result is important in the theory of metastability and, 
in the context of Metropolis dynamics, is often referred to as the 
\textit{reversibility lemma}. In that framework it is simply stated 
as the probability of reaching a configuration with energy 
larger than the one of the starting point in a time exponentially 
large in the energy difference between the final and the 
initial point. In our general it is of interest to state 
a more detailed 
result on the whole tube of trajectories overcoming this  
height level fixed \textit{a priori}. 
 
To make this result quantitative, given any $x\in\cX$ and $h,\varepsilon>0$,
for any integer $n\ge1$, we consider the tube of trajectories 
\begin{equation}
\label{staz04}
\cE^{x,h}_n:=
\{(x_0,x_1,\dots)\in \cX^{\N}:
\,x_0=x
\;\;\textrm{ and }\;\;
H(x_{n-1})+\Delta(x_{n-1},x_n)\ge H(x)+h\},
\end{equation}
which is the collection of trajectories started at $x$ whose height
at step $n$ is at least equal to the value $H(x)+h$.

\begin{proposition}
\label{t:staz03}
Let $x\in\cX$ and $h>0$. For any 
$\varepsilon\in(0,h)$, 
set
\begin{equation}
\label{staz05}
\cE^{x,h}(\varepsilon)
:=\bigcup_{n=1}^{\lfloor e^{\gb(h-\varepsilon)}\rfloor}\cE^{x,h}_n.
\end{equation}
There exists $\gb_0>0$ such that 
\begin{equation}
\label{staz06}
\P_\gb[\cE^{x,h}(\varepsilon)]
\le e^{-\gb\varepsilon/2}
\end{equation}
for any $\gb>\gb_0$.
\end{proposition}
 
In words, the set $ \cE^{x,h}(\varepsilon)$ is the set of trajectories
started at $x$ and which reach the height $H(x)+h$ at a time at most equal to 
$\lfloor\exp\left(\gb(h-\varepsilon)\right)\rfloor$. 

\section{Cycle theory in the Freidlin--Wentzell setup}
\label{lecCa}
In this section we summarize some well known facts about the 
theory of cycles, which can be seen as a handy tool to study 
the phenomenon of metastability in the Freidlin--Wentzell setup. 
Indeed, in \cite{OS1} the authors developed a peculiar approach to 
cycle theory in the framework of the Metropolis dynamics, see also 
\cite{OV}. This approach was generalized in \cite{CN03} in order 
to discuss the problem of metastability in the case of reversible 
Probabilistic Cellular Automata. In the present setup however we need the more 
general theory of cycles developed in \cite{Ca}.  We showed in \cite{CNSo} that these two approaches actually coincide 
 in the particular case of the Metropolis dynamics. 

 We recall in this section some results developed by \cite{Ca}, which will turn out to be the building bricks of our approach. 


\subsection{An alternative definition of cycles}\label{cdef}


The definition of the notion of cycle given in Section \ref{s:sg} 
is based on a property of the chain 
started at a site within the cycle itself. The point of view 
developed in [\cite{OS1}, Definition 3.1] for the Metropolis case and generalized in \cite{CN03} in the framework 
 of reversible Probabilistic Cellular Automata 
is a priori rather different. The authors introduced the notion of  
\textit{energy--cycle}, which
is defined through the height level reached by paths contained
within the energy--cycle.  
\begin{definition}
\label{defOVirr}
A nonempty set $A\subset \cX$ is an energy--cycle if and only if it is either a 
singleton or it verifies the relation
\begin{equation}
\label{2+2=5}
\max_{x,y \in A} \Phi(x,y) < \Phi(A,\cX\setminus A). 
\end{equation} 
\end{definition}
 
Even if the definitions \ref{defcy} and \ref{defOVirr} were introduced independently and in quite different contexts, it
  turns out that they actually coincide. 
 More precisely, we will prove the following result (see the proof after Proposition \ref{Cacy}): 
  \begin{proposition}
\label{OVirr}
 A nonempty set $A\subset \cX$ is a cycle if and only if it is an energy--cycle. 
\end{proposition}

After proving Proposition \ref{OVirr}, we will no longer 
distinguish the notions of cycle and of energy--cycle. 

\subsection{Depth of a cycle}
\label{s:brc}
Here we introduce the key notion of \emph{depth} of a 
cycle.

In the particular case where $D$ is a cycle, 
 a relevant property is the fact that, in the large $\gb$ limit, on an exponential scale,
  neither 
  $\tau_{D^{c}}^{x}$ nor $X_{\tau_{D^{c}}}^{x}$ depend on the starting point $x \in D$.  More precisely, 
we can formulate the following strenghtening of Proposition \ref{exgen}. 

\begin{proposition}[\cite{Ca} Proposition~4.6]
\label{indcyc}
For any cycle $C\in\cC(\cX)$, $x,y\in C$, and $z\in\cX\setminus C$
\begin{equation}\label{tpssort}
\Delta_C(x,z)=\Delta_C(y,z)=:\Delta_C(z)
\;\;\textrm{ and }\;\;
\Gamma_C(x)=\Gamma_C(y)=:\Gamma(C).
\end{equation}
\end{proposition}
The quantity $\Gamma(C)$ is the \emph{depth} of the cycle $C$.

\begin{remark}
\label{remdep}
For fixed $x$, the quantity $\Gamma_D(x)$ is monotone with respect to the inclusion, namely 
for $D,D'\subset\cX$, such that $D'\subset D$, and $x\in D'$, since $\tau^x_{\cX\setminus D'}\le\tau^x_{\cX\setminus D}$, 
from \eqref{rafuE} we deduce that 
$\Gamma_{D'}(x)\le\Gamma_{D}(x)$.
From Proposition~\ref{indcyc} it follows that for any $C,C'\in\cC(\cX)$,
 $C'\subset C$ implies $\Gamma(C')\le\Gamma(C)$.
\end{remark}
        
\subsection{Cycle properties in terms of path heights}
\label{s:cph}
In the framework 
of the study of metastability, cycles have been defined 
in terms of the height attained by paths in their interior
\cite{OS1} (see also the generalization given in \cite{CN03}).
In this section we prove the equivalence between these two approaches.

 Next we recall the following result, which links the minimal height of an exit path to the quantities we introduced previously.

\begin{proposition}[\cite{Ca} Proposition~4.12]
\label{usc0}
For any cycle $C\in\cC(\cX)$ and $y\in\cX\setminus C$
\begin{displaymath}
\min_{x\in C}[H(x)+\Delta(x,y)]
=
H(C)+\Gamma(C)+\Delta_C(y),
\end{displaymath}
where we recall the notation \eqref{fondo}.
\end{proposition}
           
The subsequent natural question is about the height that a path 
can reach within a cycle. 
We thus borrow from \cite{Ca} the following result. 
          
\begin{proposition}[\cite{Ca} Proposition~4.13]
\label{usci}
For any cycle $C\in\cC(\cX)$, $x\in C$, and $y\in\cX\setminus C$, 
there exists a path
$\omega=(\omega_{1},\ldots,\omega_{n})\in\Omega_{x,y}$
such that $\omega_{i} \in C$ for $i=1,\ldots,n-1$ and
\begin{equation}
\label{usci1}
\Phi(\omega)=H(C) + \Gamma(C) + \Delta_C(y).
\end{equation} 
For any $x,y \in C$, there is a path 
$\omega=(\omega_{1},\ldots,\omega_{n})\in\Omega_{x,y}$ such that 
$\omega_{i} \in C$ for $i=1,\ldots,n$ and
\begin{equation}
\label{dent}
\Phi(\omega) 
\leq 
H(C) 
+ 
\sup\{\Gamma(\tilde C): \tilde C\in\cC(\cX),\tilde C\subset C, 
      \tilde C \neq C \} < H(C) + \Gamma(C).
\end{equation} 
\end{proposition}
 
We stress that the right hand side term of \eqref{usci1} is infinite 
unless $y\in\partial C$. 
        
In an informal way, the first part of 
Proposition \ref{usci}, together with Proposition~\ref{usc0},
states that there exists 
a path $\omega$ contained in $C$ except for its endpoint and
joining any given $x \in C$ to any
given point $y \in \partial C$ whose cost is equal to the minimal 
cost one has to pay to exit at $y$ starting from $x$. 
Furthermore, the second part can be rephrased by saying that
one can join two arbitrary points $x$ and $y$ within $C$ by 
paying an amount which is strictly less than the minimal amount the 
process has to pay to exit from $C$; indeed, using 
Remark~\ref{remdep}, the right hand 
side of \eqref{dent} can be bounded from above 
by $H(C)+\Gamma(C)$.

We stress that this last property ensures the 
existence of \emph{at least one path} contained in the cycle connecting the two 
states and of height smaller than the one that is
necessary to exit from the cycle itself. But in general, 
there could exist other paths in the cycle, connecting the same states, 
with height larger than $H(C)+\Gamma(C)$. This is a major difference with the Metropolis case, where 
 every path contained in a  cycle 
has height smaller than the one necessary to exit the cycle itself. From this point 
of view, 
the weak reversible case  is closer to the general Freidlin--Wentzel setup than 
to the Metropolis one.

Another important property is the characterization 
of the depth of a cycle in terms of the maximal height 
that has to be reached by the trajectory to exit from a cycle. 

\begin{proposition}[\cite{Ca} Proposition~4.15]
\label{t:depth}
For any cycle $C\in\cC(\cX)$
\begin{equation}
\label{depth}
\Gamma(C) 
= 
\max_{x\in C}\big[\min_{y\in\cX\setminus C} \Phi(x,y) - H(x)\big]. 
\end{equation}  
\end{proposition}

We state now a result in which we give a different 
interpretation 
of the depth of a cycle in terms of 
the minimal height necessary to 
exit the cycle. 

\begin{proposition}
\label{fin00}
Let $C\in\cC(\cX)$  be a cycle. Then
\begin{displaymath}
\Gamma(C)=\Phi(C,\cX\setminus C)-H(C).
\end{displaymath}
\end{proposition}

\medskip
\par\noindent
\textit{Proof.\/}
Since any path connecting $C$ to $\cX\setminus C$
has at least one direct jump from a state in $C$ to a state outside of $C$, 
we have that 
\begin{displaymath}
\Phi(C,\cX\setminus C)
\ge
\min_{y\in\cX\setminus C}
\min_{x\in C}
[H(x)+\Delta(x,y)].
\end{displaymath}
Now, recalling that the principal boundary $\cB(C)$ is nonempty, by Proposition~\ref{usc0}
we have
\begin{displaymath}
\Phi(C,\cX\setminus C)\ge H(C)+\Gamma(C).
\end{displaymath}

To get the opposite bound we pick $\bar{x}\in C$ and 
$\bar{y}\in\cX\setminus C$ such that $\bar{y}\in\cB(C)$. 
Then, by the first part of Proposition~\ref{usci} there 
exists a path $\omega\in\Omega_{\bar{x},\bar{y}}$ such that 
$\Phi(\omega)=H(C)+\Gamma(C)$. 
Hence, we have that $\Phi(\bar{x},\bar{y})\le\Phi(\omega)=H(C)+\Gamma(C)$.
Finally, 
\begin{displaymath}
\Phi(C,\cX\setminus C)
=
\min_{x\in C}
\min_{y\in\cX\setminus C}
\Phi(x,y)
\le
\Phi(\bar{x},\bar{y})
\le
H(C)+\Gamma(C),
\end{displaymath}
which completes the proof.
\qed

We are now ready to discuss the equivalence between the probabilistic 
\cite{Ca} and energy \cite{OS1} approaches to cycle theory.
For any $\gl \in \R$, consider the equivalence relation
\begin{displaymath}
\cR_{\gl}:=\big\{(x,y)\in\cX^{2}, x \neq y, \Phi(x,y) < \gl \big\} 
           \cup 
           \big\{(x,x), x \in \cX \big\}. 
\end{displaymath}  
    
\begin{proposition}[ \cite{Ca} Proposition~4.18] 
\label{Cacy}
For any $\lambda\in\R$ the equivalence classes in $\cX/\cR_\lambda$ 
are either singletons $\{x\}$ such that $H(x)\ge\lambda$ or 
cycles $C\subset \cC(\cX)$ such that 
\begin{equation}
\label{Cacye0}
\max\{H(\tilde{C})+\Gamma(\tilde{C}),\, 
          \tilde C\in\cC(\cX),\,\tilde C\subset C, 
      \tilde C \neq C \}
<
\lambda
\le
H(C)+\Gamma(C).
\end{equation}
Thus we have 
\begin{equation}
\label{Cacye}
\cC(\cX)
=
\bigcup_{\lambda\in\R}\cX/\cR_\lambda.
\end{equation}
\end{proposition}

    
The results we have listed above allow us to finally prove the  equivalence between the probabilistic \cite{Ca} and energy 
approaches \cite{OS1,OV,CN03} to cycle theory, that is Proposition \ref{OVirr}. 

\medskip
\par\noindent
\textit{Proof of Proposition~\ref{OVirr}.\/}
The case $A$ is a singleton is trivial. We assume $A$ is not a singleton 
and prove the two implications. 

First assume $A$ satisfies \eqref{2+2=5}, then $A$ is an equivalence 
class in $\cX/\cR_{\Phi(A,\cX\setminus A)}$.
Thus, by Proposition~\ref{Cacy}, it follows that $A$ is a cycle. 

Reciprocally, assume that $A$ is a cycle. By \eqref{Cacye},
there exists $\lambda$ such that $A$ is an equivalence class 
of $\cX/\cR_\lambda$. Moreover, by 
\eqref{Cacye0} we have that 
\begin{displaymath}
\lambda
\le 
H(A)+\Gamma(A)
=
\Phi(A,\cX\setminus A)
\end{displaymath}
where in the last step we made use of
Proposition~\ref{fin00}.
\qed
%
 
We stress that the following 
properties are trivial in the Metropolis and in the 
weak reversible setups mentioned in Remark~\ref{partcas}, 
whereas in the general Freidlin--Wentzell setup, they are  
consequences of the non--trivial 
properties discussed previously in this section (see also \cite{CNSo}). 

For example item~\ref{fin00i02} in the following proposition
states that the principal boundary of a non--trivial cycle 
is the collection of the sites outside the cycle that can be reached 
from the interior via a single jump at height equal to the minimal 
height that has to be bypassed to exit from the cycle. 
This is precisely the notion of principal boundary adopted in 
\cite{CN03,CNS08j} in the context of reversible Probabilistic 
Cellular Automata. Note also
that such a notion 
is an obvious generalization of the idea of set of minima of the Hamiltonian
of the boundary of a cycle used in the context of Metropolis 
systems.

\begin{proposition}
\label{t:fin01}
Let $C\in\cC(\cX)$  be a cycle. Then
\begin{enumerate}
\item\label{fin00i02}
$\cB(C)=\{y\in\cX\setminus C,\,
          {\displaystyle \min_{x\in C}}
             [H(x)+\Delta(x,y)]=\Phi(C,\cX\setminus C)\}$;
\item\label{fin00i03}
$V_x<\Gamma(C)$ for any $x\in C\setminus\cF(C)$;
\item\label{fin00i04}
$V_x\ge\Gamma(C)$ for any $x\in\cF(C)$.
\end{enumerate}
\end{proposition}

\medskip
\par\noindent
\textit{Proof.\/}
Item~\ref{fin00i02}.\ This result is an immediate consequence of
Propositions~\ref{fin00} and \ref{usc0}.

Item~\ref{fin00i03}.\ 
Pick $x\in C\setminus\cF(C)$ and $y\in\cF(C)$. 
By Proposition~\ref{OVirr} we have that 
$\Phi(x,y)<\Phi(C,\cX\setminus C)$.
Thus:
\begin{displaymath}
\Phi(x,y)-H(x)<\Phi(C,\cX\setminus C)-H(x)
<\Phi(C,\cX\setminus C)-H(C),
\end{displaymath}
where we used $H(C)<H(x)$.

Item~\ref{fin00i04}.\ 
Pick $x\in\cF(C)$. 
Since $\cI_x\subset\cX\setminus C$, we have that 
$\Phi(x,\cI_x)\ge\Phi(C,\cX\setminus C)$. Since $H(x)=H(C)$, 
this entails
\begin{displaymath}
\Phi(x,\cI_x)-H(x)\ge\Phi(C,\cX\setminus C)-H(C).
\end{displaymath}
 
The item finally follows from Proposition~\ref{fin00} and 
definition \eqref{sl001}.
\qed
          
\subsection{Exit times of cycles}
\label{cet}
The main reason for which the notion of cycles has been 
introduced in the literature is that one has good control on
their exit times in the large deviation regime.
We summarize these properties in the following proposition.
    
\begin{proposition}
\label{excy}
For any cycle $C\in\cC(\cX)$, $x\in C$, 
and any 
$\gep>0$, we have that
\begin{enumerate}
\item 
\label{excyi01}
the function 
\begin{equation}\label{dalalt}
\gb \in \R^{+} \mapsto 
   \P_\beta\big[\tau_{\partial C}^{x}>e^{\gb(\Gamma(C)+\gep)}\big] 
    \end{equation}
     is SES;
\item 
\label{excyi02}
the following inequality holds for any $\delta > 0$:  
\begin{equation}
\label{dalbas}
\lim_{\beta\to\infty}
-\frac{1}{\gb} 
  \log\P_\beta\big[ \tau_{\partial C}^{x} < e^{\gb(\Gamma(C) - \delta)}\big]
\geq  \gep;
  \end{equation} 
\item 
\label{excyi03}
for any $z \in C$
\begin{equation}
\label{vis}
\lim_{\gb \to \infty} -\frac{1}{\gb} 
                \log\P_{\gb}\big[\tau^{x}_{z} > \tau^{x}_{\partial C}\big] 
 > 0;
\end{equation} 
\item 
\label{excyi04}
for any $y\in\partial C$
\begin{equation}
\label{rafucy}
\lim_{\gb \to \infty} -\frac{1}{\gb} 
      \log\P_{\gb}\big[ X_{\tau^x_{\partial C}} = y\big] 
      =\min_{x\in C}[H(x)+\Delta(x,y)] -[H(C) + \Gamma(C)].
\end{equation} 
\end{enumerate}
\end{proposition}
    
This result is the refinement of 
Proposition~\ref{exgen} in the sense that the control on the exit times and
exit locations in \eqref{rafucy}  
holds
\textit{independently} of the starting point of the process inside the cycle. 

The results of Proposition~\ref{excy} are proven in \cite{Ca}. 
More precisely, item~\ref{excyi01} is the content of 
the first part of \cite[Proposition~4.19]{Ca}. 
Item~\ref{excyi02} is \cite[Proposition~4.20]{Ca}.
Item~\ref{excyi03} is nothing but the property defining 
the cycles, see Definition~\ref{defcy} above.
Item~\ref{excyi04} follows immediately by Propositions~\ref{exgen}, 
\ref{indcyc}, and \ref{usc0}.

By combining Proposition~\ref{usc0} and 
equations \eqref{dalalt} and \eqref{rafucy} we can deduce in a trivial 
way\footnote{To deduce the corollary we use the following elementary remark: 
given two events $A,B$ such that
$(1/\gb) \log\P_\gb(B) \to 0$ and 
$(1/\gb) \log\P_\gb(A) \to - \infty$  in the limit $\gb\to\infty$,
we have that $(1/\gb)\log\P_\gb(A^{c} \cap B) \to 0$, 
where $A^{c}$ denotes the event complementary to $A$.
Indeed, since 
$\P_\gb(A^{c} \cap B) \geq \P_\gb(B) - \P_\gb(A)$, 
we get that 
$\log\P_\gb(A^{c} \cap B) \geq 
\log\P_\gb(B)
+
\log(1-\P_\gb(A)/\P_\gb(B))$. 
Then $(1/\gb)\log\P_\gb(B)\geq -\gep$ as soon as 
$\gb$ is large enough, and on the other hand, since
$\log(\P_\gb(A)/\P_\gb(B))\to-\infty$
as $\gb\to\infty$, we get that 
$\P_\gb(A)/\P_\gb(B) \to 0$.}
the following useful corollary.

\begin{corollary}
\label{corex}
For any cycle $C\in\cC(\cX)$, $\gep > 0$,
$x\in C$, and 
$y \in \cB(C)$, we have that 
          \begin{equation}\label{triv}
\lim_{\beta\to\infty}
\frac{1}{\gb} 
\log 
\P_\beta\big[\tau_{\partial C}^{x} < e^{\gb(\Gamma(C) + \gep)}, 
     X_{\tau^x_{\cX\setminus C}}= y\big] 
=0.
\end{equation} 
\end{corollary}

We discuss an interesting consequence of 
Proposition~\ref{t:staz03}.
For a given cycle $C$,
starting from the bottom of $C$, the probability of reaching an energy 
level higher
than the minimal cost necessary to exit $C$ before exiting $C$ is 
exponentially small in $\gb$. In an informal
way, this means that at the level 
of the typical behavior of trajectories, at least for trajectories 
started from $\cF(C)$, the classical 
notion of cycle for the Metropolis dynamics (which is defined in 
terms of energies only, see for example 
\cite[Chapter 6]{OV}) and the one of energy cycles are close even in 
the Freidlin--Wentzell setup. More precisely we state 
the following proposition. 
 
\begin{proposition}\label{t:quasimet}
  For any $C \in \cC(\cX)$, any $\gep > 0$ and for $\gb$ large enough: 
\begin{equation}\label{quasimet}
  \sup_{z \in \cF(C)} \P_{\gb}[ \Phi((X^{z}_{t})_{0 \leq t \leq \tau^{z}_{\cX \setminus C}}) > H(C) + \Gamma(C) + \gep ] \leq e^{-\gb \gep}.
\end{equation}      
\end{proposition}

Let us remark that we expect 
\ref{t:quasimet} to hold as well starting from anywhere within
$C$, but the proof of this result 
should be more involved.  

\subsection{Downhill or via typical jumps connected systems of cycles}
\label{tjsc}
Beside the estimate on the typical time needed to exit from a cycle, 
an important property is the one stated in \eqref{rafucy} 
which implies that when the chain exits a cycle it will 
pass typically through the principal boundary. 
This leads us to introduce the collections of pairwise disjoint cycles such that 
it is possible to go from any of them to any other 
by always performing exits through the principal boundaries.
To make this idea precise we introduce the following notion of 
oriented connection. 

\begin{definition}
\label{down}
Given two disjoint cycles 
$C,C'\in\cC(\cX)$, we say that 
$C$ is \emph{downhill connected} or \emph{ connected via typical jumps (vtj)}  to $C'$ 
if and only if $\cB(C)\cap C'\neq\emptyset$. 
\end{definition}

The fact that we introduced two names for the same notion
deserves a comment: in \cite{MNOS} 
\emph{downhill} connection is introduced in the framework of the Metropolis 
dynamics. In our opinion its natural extension to the 
general rare transition setup is the \emph{typical jumps}
connection defined in \cite[Proposition~4.10]{Ca}. 
This is the reason for the 
double name, nevertheless, in the sequel, we will always use the 
second one, which appears to be more appropriate in our setup, and we will use the abbreviation \emph{vtj}.

A \emph{vtj--connected path of cycles} 
is a 
pairwise disjoint sequence of cycles $C_1,\dots,C_n\in\cC(\cX)$ 
such that 
$C_{i}$ is 
vtj--connected to $C_{i+1}$ for all 
$i=1,\ldots,n-1$.
A \emph{vtj--connected system of cycles} 
is a 
pairwise disjoint collection of cycles $\{C_1,\dots,C_n\}\subset\cC(\cX)$ 
such that for any $1\le i<i'\le n$ there exists 
$i_1,\dots, i_m\in\{1,\ldots,n\}$ 
such that $i_1=i$, $i_m=i'$, and 
$C_{i_1},\ldots,C_{i_{m}}$ is a vtj--connected path of cycles.

We let an \emph{isolated vtj--connected system of cycles} 
to be a 
vtj--connected system of cycles
$\{C_1,\dots,C_n\}\subset\cC(\cX)$ 
such that 
\begin{displaymath}
\cB(C_i)\subset\bigcup_{j=1}^nC_j
\end{displaymath}
for any $1\le i\le n$.

Via typical jumps connected systems satisfy the following important property: 
the height that has to be reached to exit from any of the 
cycles within the system is the same. Moreover, if the system is isolated,
then the union of the 
cycles in the system is a cycle. More precisely we state 
the following two propositions. 

\begin{proposition}
\label{t:tjcs}
Let $\{C_1,\dots,C_n\}$ be a vtj--connected 
system of cycles. Then, for any $1\le i<i'\le n$, we have 
that $\Phi(C_i,\cX\setminus C_i)=\Phi(C_{i'},\cX\setminus C_{i'})$.
\end{proposition}

\medskip
\par\noindent
\textit{Proof.\/}
Consider $C_{i}$ and $C_{j}$, $1 \leq i < j \leq n$. 
By definition of a vtj--connected system, 
there exists a path of cycles consisting of 
 vtj--connected elements joining $C_{i}$ to $C_{j}$, that is 
\begin{displaymath}
C_{i} = C_{i_1},C_{i_2},\ldots,C_{i_{m-1}},C_{i_m}=C_{j}
\;\;\textrm{ such that }\;\;
\cB(C_{i_k})\cap C_{j_{k+1}}\neq\emptyset
\;\;\textrm{ for }\;\;
k=1,\ldots,m-1,
\end{displaymath}
 where all the indexes $k_{j},$ for $ j \leq i_{m},$ belong to $[1,\ldots,n]$. 
  
Now, given $k\in\{1,\ldots,m-1\}$ consider $x\in C_{i_k}$
and $y\in\cB(C_{i_k})\cap C_{i_{k+1}}$.
By Proposition~\ref{OVirr} and item~\ref{fin00i02} in 
Proposition~\ref{t:fin01} we have that 
$\Phi(x,y)=\Phi(C_{i_k},\cX\setminus C_{i_k})$.
If 
$\Phi(C_{i_{k+1}},\cX\setminus C_{i_{k+1}})>
 \Phi(C_{i_{k}},\cX\setminus C_{i_{k}})$, then 
we would have $\Phi(y,x)>\Phi(x,y)$, which is absurd in view of 
Proposition~\ref{sym}. Thus 
\begin{displaymath}
\Phi(C_{i_{k}},\cX\setminus C_{i_{k}})\ge
 \Phi(C_{i_{k+1}},\cX\setminus C_{i_{k+1}})
\end{displaymath}
for any $k=1,\ldots,m-1$. 

Iterating this inequality along the cycle 
path $\left(C_{i_1},C_{i_2},\ldots,C_{i_{m-1}},C_{i_m}\right)$, we get that 
  $\Phi(C_{i},\cX\setminus C_{i}) \geq \Phi(C_{j},\cX\setminus C_{j})$, and by symmetry we get
   \begin{equation}
    \Phi(C_{i},\cX\setminus C_{i}) \geq \Phi(C_{j},\cX\setminus C_{j}).
   \end{equation} 
Since $i$ and $j$ were chosen arbitrarily in our vtj--connected 
system, we are done. 
\qed

\begin{proposition}
\label{t:tjcs02}
Let $\{C_1,\dots,C_n\}$ be a vtj--connected 
system of cycles. 
Assume that the system is isolated
(recall the definition given above).
Then $\bigcup_{j=1}^nC_j$ is a cycle.
\end{proposition}

\medskip
\par\noindent
\textit{Proof.\/}
Let $C=\bigcup_{j=1}^nC_j$. 
From Proposition~\ref{t:tjcs}, there exists $\lambda\in\R$ such that 
$\lambda=\Phi(C_j,\cX\setminus C_j)$ for any $j=1,\ldots,n$. 

Consider $x,x'\in C$ and let $i,i'\in\{1,\ldots,n\}$ such that 
$x\in C_{i}$ and $x'\in C_{i'}$.
If $i=i'$, then by Proposition~\ref{OVirr} 
we have that $\Phi(x,x')<\lambda$. 
If, on the other hand, $i\neq i'$, by definition of 
vtj--connected system there exists 
$i_1,\dots,i_m$ such that $C_{i_k}$ is vtj--connected 
to $C_{i_{k+1}}$ for any $k=1,\ldots,m-1$.   
Thus, by using 
Proposition~\ref{OVirr} and item~\ref{fin00i02} of 
Proposition~\ref{t:fin01}, we can prove 
that $\Phi(x,x')=\lambda$. 
In conclusion, we have proven that 
$\Phi(x,x')\le\lambda$ for any $x,x'\in C$. 

Finally, since the system is isolated we have that 
$\Phi(C_i,\cX\setminus C)>\lambda$ for any $i=1,\ldots,n$
and hence, 
$\Phi(C,\cX\setminus C)>\Phi(x,x')$ for any 
$x,x'\in C$.
Thus, by Proposition~\ref{OVirr}, we have that 
$C$ is a cycle.  
\qed

\subsection{Partitioning a domain into maximal cycles}
\label{mcp}
 
In the proof of our main results a fundamental tool will be the partitioning of a set into maximal subcycles. By maximal 
we mean that given such a partition into cycles, the union of any of them 
is either the whole set or is not a cycle. 

More precisely, 
consider $D \subset \cX$ nonempty.
A \emph{partition into cycles 
of} $D$ is a partition $\{C_i,\,i\in I\}$ of $D$, 
where $I$ is a finite set of indexes, such that $C_i\in\cC(\cX)$ for any 
$i\in I$. 

\begin{definition}
\label{max}
A partition into \emph{maximal cycles} 
of the nonempty set $D\subset\cX$ is a partition 
$\{C_i,\,i\in I\}$ of strict subcycles of  $D$ 
such that the union of a number strictly smaller than $|I|$
of the cycles $C_i$'s is not a cycle. 
\end{definition}

The existence of such a partition is ensured by Proposition~\ref{tree} 
and by the fact that singletons are themselves cycles. 
In Section~\ref{exmp} 
we describe a constructive way to get such a partition for any set $D$. 

In the case where $D$ is itself a cycle, this partition 
into maximal cycles is reduced to the set $D$. 
In such a case, we can 
nevertheless decompose it into maximal \textit{strict} subcycles. 

%

\begin{proposition}[\cite{Ca} Proposition~4.10]\label{decmaxcy}
Consider a non trivial cycle $C \in \cC(\cX)$ (in particular $|C| \geq 2$), 
and its decomposition into maximal 
strict subcycles
$C = \bigsqcup_{j=1}^{n_{0}} C_{j}$ where $C_{j}$ are disjoint elements 
of $\cC(\cX)$, $n_{0} \geq 2$. 
The existence of such a decomposition is ensured by the tree 
structure of Proposition~\ref{tree}. 

The collection $\{C_1,\dots,C_{n_0}\}$ is an 
isolated 
vtj--connected system of cycles. 
Finally, from Propositions~\ref{t:tjcs} and \ref{fin00} it 
follows that 
\begin{equation}\label{constsubcy}
H(C_{i}) + \Gamma(C_{i}) = H(C_{j}) + \Gamma(C_{j})
\end{equation} 
for any $i,j \leq n_{0}$.
\end{proposition}
 
  \begin{remark}
   We stress that the original Proposition 4.10 in \cite{Ca} is actually much more exhaustive than the version presented here,
  and it allows in particular 
   to construct the set of cycles $\cC(\cX)$ in a recursive way by computing at the same time the quantities
   $\Gamma(C)$ and the $\Delta_{C}(y)$ (for $y \in \partial C$) for any
    element $C \in \cC(\cX)$, but this version will be enough 
   for our purposes. We refer to \cite{Ca} for more details. 
  \end{remark}
  
   \begin{remark}\label{selinC} For $x,y \in \cX$, from Proposition
   \ref{decmaxcy} and from Definition \ref{defgateFW}, one trivially gets the inclusion 
    \begin{displaymath}
     \cS(x,y) \subset C_{x,y}. 
    \end{displaymath} 
     
   \end{remark}

A useful property of a partition of a domain 
into maximal cycles is contained in the following proposition.
          
\begin{proposition}
\label{grtyp}
Consider a partition $\{C_i, i\in I\}$ into maximal
cycles of a nonempty set $D\subset\cX$.
Let $J\subset I$ such that 
$\{C_j,\,j\in J\}$ is a vtj-connected system of cycles. 
Then this system is not isolated, 
namely, 
there exists $j\in J$ such that 
\begin{displaymath}
\cB(C_j)\cap\Big[(\cX\setminus D)\cup\bigcup_{j'\in I\setminus J}C_{j'}
            \Big]\neq\emptyset.
\end{displaymath}
\end{proposition}

\medskip
\par\noindent
\textit{Proof.\/}
The proposition follows immediately by the maximality assumption on
the partition of $D$ and by Proposition~\ref{t:tjcs02}.
\qed

As a consequence of the above property we show that any state 
in a nonempty domain can be connected to the exterior 
of the domain by means of a vtj--connected 
cycle path made of cycles belonging to the domain itself. 
This will be a crucial point in the proof of Proposition~\ref{Main}. 

\begin{proposition}
\label{dhp}
Consider a nonempty domain $D\subset\cX$. For any state $x\in D$ 
there exists a vtj--connected cycle path 
$C_1,\ldots, C_n\subset D$ with $n\ge1$ such that 
$x\in C_1$ and $\cB(C_n)\cap(\cX\setminus D)\neq\emptyset$.
\end{proposition}

\medskip
\par\noindent
\textit{Proof.\/}
If $D$ is a cycle the statement is trivial. Assume $D$ is not a cycle
and 
consider $\{C_i,\,i\in I\}$ a partition of $D$ into 
maximal cycles. Note that $|I|\ge2$.

Now, we partition $\{C_i,\,i\in I\}$ into its \emph{maximal}
vtj--connected components
$\{C^{(j)}_k,\,k\in I^{(j)}\}$, for $j$ belonging
to some set of indexes $J$.
More precisely, 
we have the following: 
\begin{enumerate}
\item
\label{uso01}
each collection
$\{C^{(j)}_k,\,k\in I^{(j)}\}$ is a vtj--connected system of cycles;
\item
$\bigcup_{j\in J} \{C^{(j)}_k,\,k\in I^{(j)}\} 
=
\{C_i,\,i\in I\}$;
\item
$C^{(j)}_k\neq C^{(j')}_{k'}$
for any $j,j'\in J$ such that 
$j\neq j'$,
any 
$k\in I^{(j)}$, and 
$k'\in I^{(j')}$.
\item
\label{uso04}
for any $j\in J$ and 
$C\in\bigcup_{j'\in J\setminus\{j\}}\{C^{(j')}_{k'},\,k'\in I^{(j')}\}$
we have that 
$\{C^{(j)}_{k},\,k\in I^{(j)}\}\cup\{C\}$ is not a 
vtj--connected system of cycles. 
\end{enumerate}

By the property \ref{uso01} above and 
by Proposition~\ref{grtyp}, if the union of the principal boundary
of the cycles of one of those components does not 
intersect the exterior of $D$, then it necessarily intersects 
one of the cycles of one of the other components. 
Otherwise stated, for any $j\in J$
\begin{equation}
\label{pp000}
\Big(\bigcup_{k\in I^{(j)}} \cB(C^{(j)}_{k})\Big)\cap(\cX\setminus D)=\emptyset
\Rightarrow
\exists j'\in J\setminus\{j\},\, k'\in I^{(j')}:\,
\Big(\bigcup_{k\in I^{(j)}} \cB(C^{(j)}_{k})\Big)\cap C^{(j')}_{k'}
\neq\emptyset. 
\end{equation}

Now, consider $x\in D$ and $j_0\in J$ such that 
$x\in \cup_{k\in I^{(j_0)}} C^{(j_0)}_k$.
We construct a 
sequence of indexes 
$j_0,j_1,\dots\in J$ by using recursively the following rule
\begin{quote}
if 
$\left(\bigcup_{k\in I^{(j_r)}} \cB(C^{(j_r)}_{k})\right)\cap(\cX\setminus D)
=\emptyset$, choose $j\in J$ such that 
there exists $k'\in I^{(j)}$ satisfying
$\left(\bigcup_{k\in I^{(j_r)}} \cB(C^{(j_r)}_{k}) \right)
\cap C^{(j)}_{k'}\neq\emptyset$
and let $j_{r+1}=j$
\end{quote}
until the if condition above is not fulfilled. 

Note that all the indexes $j_0,j_1,\dots$ are pairwise not equal, 
namely, the algorithm above does not construct loops of 
maximal vtj--connected components. 
Indeed, if there were $r$ and $r'$ such that $j_r=j_{r'}$ 
then the union of the maximal vtj--connected components 
corresponding to the indexes $j_{r},j_{r+1},\ldots,j_{r'}$ 
would be a vtj--connected system of cycles and this is absurd 
by definition of maximal connected component (see property~\ref{uso04} above). 

Thus, since the number of maximal vtj--connected components 
in which the set $\{C_i,\,i\in I\}$ is partitioned is finite, 
the recursive application of the above rule produces 
a finite sequence of indexes $j_0,j_1,\ldots,j_{r_x}$ with 
$r_x\ge0$ such that 
$\left(\bigcup_{k\in I^{(j_{r_x})}} \cB(C^{(j_{r_x})}_{k})\right)\cap(\cX\setminus D)
\neq\emptyset$.

Finally, by applying the definition of vtj--connected system of 
cycles to each component 
$\{C^{(j_r)}_k,\,k\in I^{(j_r)}\}$ for $r=0,\ldots,r_x$
we construct a vtj--connected cycle path
$C_1,\ldots,C_n\subset D$ such that 
$C_1$ is the cycle containing $x$ and belonging to the component 
$\{C^{(j_0)}_k,\,k\in I^{(j_0)}\}$ 
and 
$C_n$ is one of the cycles in the component 
$\{C^{(j_{r_x})}_k,\,k\in I^{(j_{r_x})}\}$
such that $\cB(C_n)\cap(\cX\setminus D)\neq\emptyset$.
\qed

\subsection{Example of partition into maximal cycles}
\label{exmp}
It is interesting 
to discuss a constructive way to exhibit a partition into maximal cycles of a given $D \subset \cX$. 
For this reason we now describe a method 
inherited from the Metropolis setup in \cite{MNOS}.
For $D \subset\cX$ nonempty and $x\in D$, we consider
\begin{equation}
\label{excR}
R_{D}(x): =\{x\}\cup\{y\in\cX, \Phi(x,y) < \Phi(x,\cX\setminus D) \},
\end{equation}   
namely, $R_D(x)$ is the union of $\{x\}$ and of the points in $\cX$ 
which  
can be reached by means of paths starting from $x$ with height smaller that the 
height that it is necessary to reach  
to exit from $D$ starting from $x$. 
    
\begin{proposition}
\label{stablvl}
Given the nonempty set $D\subset\cX$ and $x\in D$, 
\begin{enumerate}
\item  
the following inclusion holds: $R_{D}(x) \subset D$;
\item 
the set $R_{D}(x)$ is a cycle;
\item 
if $x' \in  R_{D}(x) $, 
then $R_{D}(x) =  R_{D}(x')$.         
\end{enumerate}
\end{proposition}
  
\medskip
\par\noindent
\textit{Proof.\/}
The first item is clear by the definition of communication heights. 
Indeed, by contradiction, assume that there exists
$y \in R_{D}(x) \cap(\cX\setminus D)$, 
then $\Phi(x,y)$ satisfies simultaneously 
\begin{displaymath}
\Phi(x,y) \geq \Phi(x,\cX\setminus D)
\;\;\textrm{ and }\;\;
\Phi(x,y) < \Phi(x,\cX\setminus D),
\end{displaymath}
which is absurd.
    
Second item.\ We consider 
$u,v \in R_{D}(x)$ and we show that 
$\Phi(u,v) < \Phi(x,\cX\setminus A)$. As a consequence, we will get  
that $R_{D}(x)$ is a maximal connected subset of 
$\cX$ satisfying that the maximum internal communication cost is 
strictly smaller than the given threshold $\Phi(x,\cX\setminus D)$, 
and, by Proposition~\ref{Cacy}, these sets are cycles. 

We use a concatenation argument. Namely, consider
$\go\in\Omega^\textrm{opt}_{u,x}$ and 
$\go'\in\Omega^\textrm{opt}_{x,v}$ and let 
$\go'' \in\Omega_{u,v}$ be 
the path obtained by concatenating $\go$ and $\go'$.
We then have 
\begin{displaymath}
\Phi(u,v) 
\leq \Phi(\go) \vee \Phi(\go')
\end{displaymath} 
and hence 
\begin{displaymath}
\Phi(u,v) \leq \Phi(u,x) \vee \Phi(x,v).
\end{displaymath} 
 
By the symmetry property in Proposition~\ref{sym}, 
we get that $\Phi(u,x) = \Phi(x,u)$.
Since by construction 
$\Phi(x,u) < \Phi(x,\cX\setminus D)$ and 
$\Phi(x,v) < \Phi(x,\cX\setminus D)$, we 
get indeed $\Phi(u,v) < \Phi(x,\cX\setminus D)$.
     
Third item.\ 
We first claim that 
\begin{equation}\label{eqH}
\Phi(x',\cX\setminus D)
=
\Phi(x,\cX\setminus D) 
\;\;\;\textrm{ for any } x' \in R_{D}(x).
    \end{equation} 
     
To prove \eqref{eqH} pick $x'\in R_D(x)$.
First assume that 
$\Phi(x',\cX\setminus D) < \Phi(x,\cX\setminus D)$. 
Then, 
we can consider 
a path  
$\go\in\Omega_{x,x'}$ 
such that $\Phi(\go) < \Phi(x,\cX\setminus D)$ 
and a path $\go'\in\Omega^\textrm{opt}_{x',\cX\setminus D}$. 
Note that 
$\Phi(\omega')=\Phi(x',\cX\setminus D)<\Phi(x,\cX\setminus D)$.
Now, by concatenation of the two preceding paths, we obtain a 
path $\go''\in\Omega_{x,\cX\setminus D}$ such that 
$\Phi(\go'') = \Phi(\go) \vee \Phi(\go') < \Phi(x,\cX\setminus D)$, 
which is absurd. Hence, we have that 
$\Phi(x',\cX\setminus D) \geq \Phi(x,\cX\setminus D)$.
       
To prove the opposite inequality, consider 
$\go \in\Omega^\textrm{opt}_{x',x}$. 
From the Proposition~\eqref{sym}, 
we get that 
$\Phi(\go)=\Phi(x',x)=\Phi(x,x')<\Phi(x,\cX\setminus D)$. 
Similarly, consider a path 
$\go' \in\Omega^\textrm{opt}_{x,\cX\setminus D}$ and note that 
$\Phi(\go')=\Phi(x,\cX\setminus D)$. 
Then, the path $\go''\in\Omega_{x',\cX\setminus D}$ 
obtained by concatenating $\go$ and $\go'$ satisfies 
$\Phi(\go'')=\Phi(x,\cX\setminus D)$, from which we deduce 
$\Phi(x',\cX\setminus D) \leq \Phi(x,\cX\setminus D)$. 
The proof \eqref{eqH} is thus completed. 
    
Now we come back to the proof of the third item.
We consider $x'\in R_{D}(x)$
and proceed by double inclusion.
We first show that $R_{D}(x')\subset R_{D}(x)$. 
Pick up $y \in R_{D}(x')$: 
from the definition of $R_{D}(x')$ and \eqref{eqH}, we get that 
$\Phi(x',y) <\Phi(x',\cX\setminus D)= \Phi(x,\cX\setminus D)$. 
Now we consider 
$\go\in\Omega^\textrm{opt}_{x,x'}$, and by a concatenation 
argument similar to the one we already used twice, 
we get that 
\begin{displaymath}
\Phi(x,y)\leq\Phi(\go)\vee\Phi(x',y)<
\Phi(x,\cX\setminus D)\vee\Phi(x',\cX\setminus D)=\Phi(x,\cX\setminus D),
\end{displaymath}
which implies $R_{D}(x')\subset R_{D}(x)$. 

On the other hand, the inclusion $R_{D}(x)\subset R_{D}(x')$ proceeds in 
the same vein.
Consider $y\in R_{D}(x)$ so that $\Phi(x,y) < \Phi(x,\cX\setminus D)$.
Pick up a path $\go\in\Omega^\textrm{opt}_{x',x}$. 
Using again the symmetry of $\Phi$, 
we get that $\Phi(\go)=\Phi(x',x)=\Phi(x,x')<\Phi(x,\cX\setminus D)$.
Moreover, a concatenation argument shows that 
\begin{displaymath}
\Phi(x',y) 
\leq \Phi(\omega) \vee \Phi(x,y) 
< \Phi(x,\cX\setminus D) 
\end{displaymath} 
where we have also used that $y\in R_{D}(x)$.
Finally, 
from \eqref{eqH},
we deduce 
$\Phi(x',y) < \Phi(x',\cX\setminus D)$, which implies $y\in R_{D}(x')$.
\qed
    
The main motivation for introducing the sets 
\eqref{excR} is the fact that they provide in a constructive way 
a partition of a given set into maximal subcycles. 
The existence of such a partition is ensured by the 
structure of the set of cycles, see Proposition~\ref{tree}, but 
we point out that this way of obtaining the maximal subcycles of 
a given set $D$ seems to be new in the context of the irreversible dynamics. 
Before stating precisely this result, for $D\subset\cX$, we set 
\begin{equation}
\label{setcr}
\cR_D:=
\{C\in\cC(\cX),\,
  \textrm{ there exists } x\in D\textrm{ such that } C=R_{D}(x)\}.
\end{equation}

\begin{proposition}
\label{maxdec}
Let $D \subset \cX$ nonempty, then $\cR_D$ is a 
partition into maximal cycles of $D$.
\end{proposition}
   
\medskip        
\par\noindent
\textit{Proof.\/}
In view of definition \eqref{excR} and Proposition~\ref{stablvl}, 
the only not obvious point of this result is the one concerning 
maximality. 
Note that the maximality condition on cycles can be stated
equivalently as follows: any cycle $C\in\cC(\cX)$ such that there 
exists $R\in\cR_D$ verifying $R\subset C$ and $R\neq C$ satisfies 
$C \cap(\cX\setminus D)\neq \emptyset$.

Now, 
assume that $C\in\cC(\cX)$ is a cycle strictly containing 
$R_{D}(x)$ for some $x\in D$. 
We will show that necessarily $C\cap(\cX\setminus D)\neq\emptyset$. 

By definition of $R_{D}(x)$, $C$ contains 
a point $v \notin R_{D}(x)$, that is $\Phi(x,v) \geq \Phi(x,\cX\setminus D)$. 
As both $x$ and $v$ are elements of $C$, recalling Proposition~\ref{OVirr},
we get that 
\begin{displaymath}
\Phi(C,\cX\setminus C) > \Phi(x,v) \geq \Phi(x,\cX\setminus D).
\end{displaymath}
 
On the other hand, we can choose $y \in\cX\setminus D$ 
such that there exists $\go\in\Omega_{x,y}$ satisfying
$\Phi(\go) = \Phi(x,\cX\setminus D)$. 
Then the above bound implies that 
$\Phi(C,\cX\setminus C)>\Phi(\go)$ 
and 
in particular $y\in\cX\setminus D$. Hence the result. 
\qed

\section{Proof of main results}
\label{s:pmr}
In this section we prove the results stated in 
Sections~\ref{s:meta} and \ref{s:colla}. 
We follow the scheme of\cite{MNOS}, 
but the proofs are a bit different. 
The proofs of Theorems~\ref{tune} and \ref{expbeh} 
 are  quite
similar to the analogous ones in \cite{MNOS}, nevertheless we chose to include them for the sake of completeness. 

\medskip
\par\noindent
\textit{Proof of Theorem~\ref{t:cp}.\/}
Proof of \eqref{cp01}.\ 
Let $C$ be the set of states $y\in\cX$ such that 
$\Phi(x,y)<V^\textrm{m}+H(x)$. 
By Proposition~\ref{OVirr} the set $C$ is a cycle and, 
by construction, $x\in\cF(C)$ and 
$\Phi(C,\cX\setminus C)=V^\textrm{m}$.
Hence, 
by Proposition~\ref{fin00} it follows that $\Gamma(C)=V^\textrm{m}-H(x)$.
Finally, since $\cX^\textrm{s}\cap C=\emptyset$ 
implies $\tau^x_{\cX^\textrm{s}}\ge\tau^x_{\partial C}$, 
we have that \eqref{cp01} 
follows by item~\ref{excyi02} in Proposition~\ref{excy}.

Proof of \eqref{cp02}.\ 
As we have already remarked at the end of Section~\ref{defmet}, see also 
Figure~\ref{f:mestru}, $\cX_{V^\textrm{m}}=\cX^\textrm{s}$. 
Hence, \eqref{cp02} is an immediate consequence of 
Proposition~\ref{Main}.
\qed
      
\medskip
Before proving Theorem~\ref{tune} we first state and prove the 
following preliminary integrability result. 

\begin{lemma}\label{uihit}
Given any real $\gd > 0$ and any state $x \in \cX$, the family of random 
variables 
$Y_{\gb}^{x} = 
 \tau_{\cX^{\normalfont\textrm{s}}}^{x} 
  e^{-\gb(V^{\normalfont\textrm{m}}+\gd)}$ is 
uniformly integrable, more precisely,
for any $n\ge1$
\begin{equation}
\label{e:uihit}
\sup_{x \in \cX} \P_\gb[\tau_{\cX^{\normalfont\textrm{s}}}^{x} 
      e^{-\gb(V^{\normalfont\textrm{m}}+\gd)} > n ] 
\le
\frac{1}{2^n}
\end{equation}
for $\beta$ large enough.
\end{lemma}
  
\medskip
\par\noindent
\textit{Proof.\/}
For any $n \geq 1$, by making use of the Markov property, we directly get
\begin{displaymath}
\sup_{x \in \cX} \P_\gb[\tau_{\cX^{\textrm{s}}}^{x} 
      e^{-\gb(V^{\textrm{m}}+\gd)} > n ] 
\leq 
\left( \sup_{x' \notin \cX^{\textrm{s}}} 
  \P_\gb[\tau_{\cX^{\textrm{s}}}^{x}  > e^{\gb(V^{m}+\gd)} ] \right)^{n}.
\end{displaymath}
Recalling that
$\cX_{V^\textrm{m}} = \cX^\textrm{s}$
(see the end of Section~\ref{defmet}) and 
making use of Proposition~\ref{Main}, we get that the 
above quantity is bounded from above 
by $2^{-n}$ as soon as $\gb$ large enough.  
\qed
           
\medskip
\par\noindent
\textit{Proof of Theorem \ref{tune}.\/}
Fix $x \in \cX^{\textrm{m}}$ and $\gd > 0$.
Combining the convergence to zero in probability of the random variables 
$Y_{\gb} = \tau_{\cX^{\textrm{s}}}^{x} e^{-(V^{\textrm{m}}+\gd)\gb}$,
which has been shown in Theorem~\ref{t:cp}
and their uniform summability stated in Lemma~\ref{uihit},
we get that the family of random variables $Y_{\gb}$ 
converges to $0$ in $L^{1}$. 
Hence, 
\begin{equation}
\label{tune000}
\E_\gb[\tau_{\cX^{\normalfont\textrm{s}}}^{x}] < e^{\gb(V^{\textrm{m}}+\gd)}
\end{equation}
for $\gb$ large enough,
          
On the other hand, by making use of the Markov's inequality we get
the following bound: 
\begin{displaymath}
\P_\gb[ \tau_{\cX^{\textrm{s}}}^{x} > e^{\gb(V^{\textrm{m}}-\gd)}] 
  \leq \E_\gb[\tau_{\cX^{\textrm{s}}}^{x}]\, e^{-\gb(V^{\textrm{m}}-\gd)}.
\end{displaymath}
Using once again Theorem~\ref{t:cp}, 
we obtain that there exists $K > 0$ such that
\begin{equation}
\label{tune010}
\E_\gb[\tau_{\cX^{\textrm{s}}}^{x}] 
  \geq e^{\gb(V^{\textrm{m}}-\gd)} \left( 1 - e^{-\gb K} \right)
\end{equation} 
as soon as $\gb$ is large enough. 
          
The Theorem~\ref{tune} finally follows from bounds \eqref{tune000} 
and \eqref{tune010}.
\qed
         
\medskip
\par\noindent
\textit{Proof of Theorem \ref{expbeh}.\/}
We first prove item 1.
Let $x_0$ be the recurrent state of Theorem \ref{expbeh}
and recall \eqref{le}--\eqref{fr}.
We consider $s,t>0$ and let 
$\tau^{x_0}_{*}(t) 
  = 
  \inf\{ n \geq tT_{\gb},\, X_{n} \in \{x_{0},\cX^{\textrm{s}} \}\}$ 
be the first hitting time to the set 
$\{x_{0},\cX^{\textrm{s}}\}$ after time $t T_{\gb}$ for the 
chain $X_n$ started at $x_0$. 
 
Then we decompose:
\begin{displaymath}
\begin{array}{rcl}
\P_\gb[ \tau_{\cX^{\textrm{s}}}^{x_{0}} > (t+s)T_{\gb}]
&\!\!=&\!\!  
\P_\gb[ \tau_{\cX^{\textrm{s}}}^{x_{0}} > (t+s)T_{\gb}; 
     \tau^{x_0}_{*}(t)  \leq t T_{\gb} + T_{\gb}'] \vphantom{\bigg\{}\\
&&\!\! 
+ \P_\gb[ \tau_{\cX^{\textrm{s}}}^{x_{0}} > (t+s)T_{\gb}; 
\tau^{x_0}_{*}(t)  > t T_{\gb} + T_{\gb}']\,\,.
\\
\end{array}
\end{displaymath}
            
Using the Markov property and of the fact that
$\{ \tau_{\cX^{\textrm{s}}}^{x_{0}} > \tau^{x_0}_{*}(t) \} 
   \subset \{X_{\tau^{x_0}_{*}(t)} = x_{0} \}$, we directly get:
\begin{equation}
\label{decomp}
\begin{array}{l}
\P_\gb[ \tau_{\cX^{\textrm{s}}}^{x_{0}} > (t+s)T_{\gb}; 
       \tau^{x_0}_{*}(t)  \leq t T_{\gb} + T_{\gb}']  
\vphantom{\bigg\{}
\\
\phantom{mmmmmmm}
{\displaystyle
    = \sum_{n=0}^{T_{\gb}'} \P_\gb[ \tau^{x_0}_{*}(t) 
 = t T_{\gb} +n,\tau_{\cX^{\textrm{s}}}^{x_{0}} > tT_{\gb} +n ] 
   \,\P_\gb[\tau_{\cX^{\textrm{s}}}^{x_{0}} > sT_{\gb} -n ].   
}
\\
\end{array}
\end{equation} 
       
Combining monotonicity and the fast recurrence property \eqref{fr}, 
by the decomposition \eqref{decomp} we deduce 
 \begin{equation}\label{lob}
  \begin{aligned}
   & \P_\gb[ \tau_{\cX^{\textrm{s}}}^{x_{0}} > (t+s)T_{\gb}; 
       \tau^{x_0}_{*}(t)  \leq t T_{\gb} + T_{\gb}']  \\
      & \phantom{m} \geq  \left( \P_\gb[ \tau_{\cX^{\textrm{s}}}^{x_{0}} 
         > t T_{\gb} + T_{\gb}'] -
      \P_\gb[ \tau_{\cX^{\textrm{s}}}^{x_{0}} > t T_{\gb} + T_{\gb}';
            \tau^{x_0}_{*}(t)  > t T_{\gb} + T_{\gb}' ]
      \right)\P_\gb \left[\tau_{\cX^{\textrm{s}}}^{x_{0}} > sT_{\gb} \right] \\
       & \phantom{m} \geq  \left( \P_\gb[ \tau_{\cX^{\textrm{s}}}^{x_{0}} 
           > t T_{\gb} + T_{\gb}'] -
      \gd_{\gb}
      \right)\P_\gb \left[\tau_{\cX^{\textrm{s}}}^{x_{0}} > sT_{\gb} \right]. 
  \end{aligned}
 \end{equation}

      Here and later, we made use of the following obvious monotonicity 
property: for $b,c \in \R$ such that $b\geq c$,    
\begin{displaymath}
\{ T \geq b \} \subset \{ T \geq c\}
\end{displaymath}      
where $T$ is any random variable. 

 We bound the same quantity from above in a similar fashion. Namely, using \eqref{decomp} once again:
\begin{equation}\label{arec}
\P_\gb[ \tau_{\cX^{\textrm{s}}}^{x_{0}} > (t+s)T_{\gb}; 
   \tau^{x_0}_{*}(t)  \leq t T_{\gb} + T_{\gb}'] \leq 
 \P_\gb \left[\tau_{\cX^{\textrm{s}}}^{x_{0}} > tT_{\gb} \right]
    \left( \P_\gb \left[\tau_{\cX^{\textrm{s}}}^{x_{0}} > sT_{\gb} - T_{\gb}' 
  \right] + \gd_{\gb}\right).
\end{equation} 
    
Consider $\gb$ large enough so that $T_{\gb}' \leq T_{\gb}$. 
For any given integer $k \geq 1$, combining
\eqref{arec} and monotonicity, we get:
\begin{displaymath}
\P_\gb[ \tau_{\cX^{\textrm{s}}}^{x_{0}} > (k+2)T_{\gb}] 
\leq \P_\gb[ \tau_{\cX^{\textrm{s}}}^{x_{0}} > kT_{\gb}] \left(\gd_{\gb} +
\P_\gb[ \tau_{\cX^{\textrm{s}}}^{x_{0}} > T_{\gb}] \right).
   \end{displaymath} 
   
Given the definition of $T_{\gb}$ (see \eqref{le}), 
there exists $r \in (0,1)$ such that 
$\gd_{\gb} + \P_\gb[ \tau_{\cX^{\textrm{s}}}^{x_{0}} > T_{\gb}] \leq r$ 
as soon as $\gb$ is large enough. As a consequence, for $\gb$ large enough, the following inequality holds:
\begin{equation}\label{uip}
\P_\gb[ \tau_{\cX^{\textrm{s}}}^{x_{0}} > kT_{\gb}] \leq r^{k/2},
\end{equation} 
and this implies the tightness of the family 
$\tau_{\cX^{\textrm{s}}}^{x_{0}}/T_{\gb}$.
   
Combining the upper bound \eqref{arec} 
and the lower bound \eqref{lob}, we deduce that the limit in law $X$ of any subsequence 
$\left(\tau_{\cX^{\textrm{s}}}^{x_{0}}/T_{\gb}\right)_{\gb_{k}}$ satisfies the relation:
\begin{equation}\label{ex}
\P_\gb[X > t+s] = \P_\gb[X > t]\,\P_\gb[X > t]
\end{equation} 
for any $t,s \geq 0$ which are continuity points for the distribution of 
$\tau_{\cX^{\textrm{s}}}^{x_{0}}$. Since the set of such points is dense in $\R$ and a distribution function is 
 always right continuous, \eqref{ex} is valid for every 
$s,t \geq 0$. This implies that $\P_\gb(X > t) = e^{-at}$ with $a \in (0,\infty]$. 
It is clear that the case $a = \infty$
is excluded from the definition of $T_{\gb}$, since it would imply 
that $X$ is
almost surely equal to zero, which is in contradiction with 
the fact that 
\begin{equation}\label{dalb}
\P_\gb[X < 1] = \lim_{\gb \to \infty} 
  \P_\gb[\tau_{\cX^{\textrm{s}}}^{x_{0}} < 1] \leq 1 -e^{-1}.
\end{equation} 
   
By the Porte--Manteau theorem, we get that 
\begin{equation}\label{dalal}
1 -e^{-1} \leq \lim_{\gb \to \infty}  
   \P_\gb[\tau_{\cX^{\textrm{s}}}^{x_{0}} \leq 1] =  \P_\gb[X \leq 1],
\end{equation} 
and combining \eqref{dalb} and \eqref{dalal}, we conclude that $a = 1$.

As for item 2, combining the dominated convergence theorem and the uniform 
summability \eqref{uip}, we can write 
\begin{displaymath}
\lim_{\gb \to \infty} \frac{\E_\gb[\tau_{\cX^{\textrm{s}}}^{x_{0}} ]}{T_{\gb}} 
    =  \lim_{\gb \to \infty} \int_{0}^{\infty}  
      \P_\gb[\tau_{\cX^{\textrm{s}}}^{x_{0}} \geq T_{\gb}t]\, \textrm{d}t 
     = \int_{0}^{\infty} \lim_{\gb \to \infty} 
      \P_\gb[\tau_{\cX^{\textrm{s}}}^{x_{0}} \geq T_{\gb}t]\,\textrm{d}t = 1,
    \end{displaymath} 
which entails the convergence \eqref{eq:expbeh}.

Item 3 directly follows from items 1 and 2 of the current theorem, which concludes the proof. 
\qed

 Now, given $x,y\in\cX$, we consider a minimal gate  $W \subset \cW(x,y)$ as in Definition \ref{defgateFW} and
 we go to the proofs of Theorem \ref{crossgate} and of Proposition \ref{traje}.

 
 To prove both these results,
 we first construct in a more formal way the tube of typical trajectories $\cK_{x,y}$ introduced in Section \ref{s:meta};  we stress that this task is performed 
 by making an
 extensive use of the notions developed in the previous parts. Then we show that $\cK_{x,y}$ is indeed typical in the low temperature
  regime, that is we show Proposition \ref{traje}.
  Our 
  task to prove Theorem \ref{crossgate} will then be reduced to show
  the inclusion  $ \cK_{x,y} \subset \{\tau^x_W<\tau^x_y\}$. 
  
   To  give an explicit description of 
 the set $\cK_{x,y}$, 
we first need to introduce some typical events and recurrent notations. 
 
We introduce the positive quantity 
\begin{displaymath}
\gd_{0} = \inf_{\go \in \Omega_{x,y}\setminus\Omega_{x,y}^\textrm{opt}}
         [\Phi(\go) - \Phi(x,y)]
        \end{displaymath} 
and let $\gep \in (0,\gd_{0}/2)$. 
We define the cycle 
\begin{displaymath}
C:=\{u \in \cX, \Phi(x,u) <  \Phi(x,y) + \gd_{0}/2 \}.
\end{displaymath}

 Of course, the cycle $C$ coincides with the cycle $C_{x,y}$ of Definition \ref{defgateFW} and 
 we define it in this way for technical purposes only.  
 
Note that 
any path $\Omega_{x,y}^{\rm{opt}}$ is contained in $C$. 
Also, we already noted (and
this is actually the major technical difference with the analogous 
result of \cite{MNOS}) that there might be paths contained in $C$, 
joining $x$ to $y$ and
which do not belong to $\Omega_{x,y}^{\rm{opt}}$.
   
We introduce the decomposition 
$\cM = \{C_{j}, j \leq n_{0} \}$ of $C$ into maximal strict 
subcycles of $C$.  The decomposition 
     $\cM$ is an isolated vtj--connected system of cycles (see Proposition \ref{decmaxcy}).
Then we 
discuss some geometrical properties of the decomposition 
$\cM$.
 
 We first note that it is clear that $x$ and $y$ are not 
contained in the same element of $\cM$. 
 Indeed, if they were contained in a common element $\bar{C} \in \cM$, we would have 
$\Phi(x,y)<H(\bar{C})+\Gamma(\bar{C})$ and in particular, 
from the definition of $C$, 
this would imply $C\subset \bar{C}$, which is absurd from 
the non triviality of the decomposition $\cM$.  
Thus we can denote by 
$C(x)$ and $C(y)$ the two (distinct) 
elements of $\cM$ containing respectively 
the states $x$ and $y$. More generally, for any $u \in C$, we define  $C(u)$ as being 
     the element of $\cM$ containing $u$.

To define $\cK_{x,y}$, we shall start to restrict the set of trajectories  
to the set of trajectories 
$\Omega_{x,y}\cap\{\tau^{x}_{y}<\tau^{x}_{\cX \setminus C}\}$, 
for which the events we are going to 
introduce are well defined. 

More precisely, 
for a given trajectory of the canonical process 
$ \omega
  \in \Omega_{x,y} \cap \{\tau^{x}_{y} < \tau^{x}_{\cX \setminus C} \}$, 
we first define 
$\theta^{x}_{0} := 0, C^{x}_{0} = C(x) $ and for $j \geq 1$:
\begin{displaymath}
\theta^{x}_{j} := \inf \left\{k \geq \theta^{x}_{j-1},\, \omega_{k} 
                      \notin C^{x}_{j-1} \right\}
\end{displaymath} 
and $C^{x}_{j} = C(\omega_{\theta^x_{j}})$ is the element of $\cM$ containing 
$\omega_{\theta^x_{j}}$. This construction goes on 
   as long as $ j \leq j_{x,y}$, where we consider 
   \begin{displaymath}
    j_{x,y} := \inf\{j \geq 1, C^{x}_{j} = C(y) \}. 
   \end{displaymath} 
   
    More generally, for any $u \in C$, we introduce the similar quantities $(\theta^{u}_{j})_{j}, (C^{u}_{j})_{j}$, with 
    notations which are self explanatory. 
    
 Then we introduce the event
\begin{displaymath}
 E_{x,y} := 
\left\{\omega\in\Omega_{x,y},\, 
\tau^{x}_{y} < \tau^{x}_{\cX \setminus C} 
\;\;\textrm{and}\;\;
\tau^{x}_{y} < \inf \left\{ k \geq \tau^{x}_{C(y)},
 \, \omega_{k} \notin C(y) \right\} \right\},
  \end{displaymath}   
which is the event that the process hits $y$ after entering $C(y)$ 
before leaving $C(y)$ for the first time.
     
For $ 0 \leq i \leq n_{0}$ and $u \in C_{i}$, we introduce the event 
\begin{displaymath}
A^{u}_{C_{i}}: = \Big\{ \omega^u_{\tau^{u}_{\cX \setminus C_{i}}} 
                        \in \cB(C_{i}) \Big\},
\end{displaymath} 
 where $(\omega^{u}_{k})_{k\geq 0}$ denotes a trajectory of the canonical process starting from $u$. 
Finally we can define the set  $\cK_{x,y}$, the tube of trajectories of 
the dynamics on its transition between $x$ and $y$: 
\begin{equation}
\label{dimo080}
\cK_{x,y} := \Big\{\omega\in E_{x,y}, \,
\bigcap_{i=0}^{j_{x,y}} A^{\omega_{\theta^{x}_{i}}}_{C^{x}_{i}} 
\Big\}.
\end{equation} 

 We refer to Section \ref{s:meta} for an informal definition of $\cK_{x.y}$. 
 

\medskip
\par\noindent
\textit{Proof of Proposition~\ref{traje}.\/}

  We prove that 
\begin{displaymath}
  \P_{\gb}[\cK_{x,y}] \geq 1 - e^{-\gb \gep}
\end{displaymath} 
  as soon as $\gb$ is large enough. 
   
  Our proof first 
relies on the fact that, given $\gd > 0$,  
    for $\gb$ large:
     \begin{equation}\label{ott1}
      \inf_{ \bar{C} \in \cM}\inf_{u \in \bar{C}} \P_{\gb}[A^{u}_{\bar{C}} ] \geq 1 - e^{- \gb \gd },
     \end{equation} 
      which follows from the finiteness of $\cX$ and Corollary \ref{corex}.
      Then we will use the fact that for any $\gep' > 0$, as soon as $\gb$ is large enough: 
     \begin{equation}\label{ott2}
      \P_{\gb}[ j_{x,y} > e^{\gep \gb}] \leq e^{- \gb \gep'},
     \end{equation} 
 which we show at the end of the proof of Proposition \ref{traje}.

   Let us note that in \cite{CaCe}, the authors showed a result related to ours, in the sense that they provide the precise 
    cost on a large deviation scale of \textit{not following a path contained in} $\cK_{x,y} \cap \Omega_{x,y}^{\rm{opt}}$ on the transition 
     from $x$ to $y$. For our sake such a level of precision is not needed. On the other hand, we had to deal with the (easy) problem 
      of giving an upper bound on the random variable $ j_{x,y}$, which was overcome in \cite{CaCe} by the notion of \textit{pruning tree}. 
      
    We show how to deduce Proposition \ref{traje} from combining \eqref{ott1} and \eqref{ott2}. 
    For lightness of notations, we introduce the conditional probability 
     \begin{displaymath}
     \widetilde{\P_{\gb}}[\cdot] := \P_{\gb}\left[\hspace{6 pt} \cdot \hspace{6 pt}
     \Big| E^{x,y} \cap \left\{ \tau^{x}_{y} < \tau^{x}_{\cX \setminus C} \right\}  \right] 
     \end{displaymath} 
   in the next 
   sequence of inequalities. Of course, since $y \in C$ and $y \in C(y)$, applying the strong Markov property at time 
     $\tau^{x}_{C(y)}$ and Definition \ref{defcy} we immediately 
    get that, for any $\gep' > 0$, as soon as $\gb$ is large enough: 
     \begin{equation}\label{condtil}
      \P_{\gb}\left[ E^{x,y} \cap \left\{ \tau^{x}_{y} < \tau^{x}_{\cX \setminus C} \right\} \right] \geq 1 - e^{-\gb \gep'}.       
     \end{equation} 
     
      It follows from this inequality that similar inequalities to \eqref{ott1} and \eqref{ott2} also hold for 
       the probability $\widetilde{\P_{\gb}} $ instead of $\P_{\gb}$, and we will still refer to these slightly modified 
       versions of \eqref{ott1} and \eqref{ott2}  as \eqref{ott1} and \eqref{ott2} in the following. 
     
      Denoting by $\gep'$
      a (small) positive constant which may change from line to line, we then get:
    
     \begin{align*}
       \P_{\gb}\left[ \cK_{x,y}\right]  
          & \geq  \widetilde{\P_{\gb}}\left[ \bigcap_{j=1}^{ j_{x,y}} A^{X^{x}_{\theta^{x}_{i}}}_{C^{x}_{i}} \right](1- e^{- \gb \gep'}) \\
            & \geq  \widetilde{\P_{\gb}}\left[ \bigcap_{j=1}^{ j_{x,y}} A^{X^{x}_{\theta^{x}_{i}}}_{C^{x}_{i}} ,  j_{x,y} \leq e^{\gb \gep} \right](1- e^{- \gb \gep'}) 
           \end{align*}
         \begin{align*}  
        \phantom{iiiiiiiiiiiiiiiiiiiiiiiiiii}  &  \geq \widetilde{\P_{\gb}}\left[  \bigcap_{j=1}^{e^{\gb \gep}} A^{X^{x}_{\theta^{x}_{i}}}_{C^{x}_{i}}, j_{x,y} \leq e^{\gb \gep}  \right]  (1 - e^{- \gb \gep'}) \\
          & \geq \left( \widetilde{\P_{\gb}}\left[  \bigcap_{j=1}^{e^{\gb \gep}} A^{X^{x}_{\theta^{x}_{i}}}_{C^{x}_{i}}  \right] - \widetilde{\P_{\gb}}[ j_{x,y} > e^{\gb \gep}]\right) (1  - e^{- \gb \gep'})  \\
          & \geq  \left( \prod_{j=1}^{e^{\gb \gep}} \inf_{\bar{C} \in \cM} \inf_{u \in A} \widetilde{\P_{\gb}}[A^{u}_{\bar{C}} ] - e^{- \gb \gep'}\right)(1  - e^{- \gb \gep'})
     \end{align*}
    where we used \eqref{condtil}, \eqref{ott2} and the strong Markov property. Now from \eqref{ott1}, we get
     \begin{displaymath}
      \begin{aligned}
      \P_{\gb}[ \tau_{y}^{x} > \tau_{W}^{x}] & \geq \left((1-e^{-\gb \gd})^{e^{\gep \gb}} - e^{- \gb \gep'} \right) (1-  e^{- \gb \gep'}) \\
          & \geq \left(e^{-e^{\gb(\gep - \gd)}} -  e^{- \gb \gep'}\right)(1-  e^{- \gb \gep'}),
      \end{aligned}
     \end{displaymath} 
     and considering $\gd > \gep$, the statement of Proposition \ref{traje} follows. 
     
    Now we are left with the proof of \eqref{ott2}.

     Since $\cM$ is an isolated vtj--connected system of cycles, we deduce that 
\begin{equation}\label{unifdentroC}
 \max_{\bar C \in \cM} \max_{u \in \bar C} \P_{\gb} \left[ C(y) \notin (C^{u}_{1}, \ldots, C^{u}_{n_{0}})\right] \leq 1 - e^{-\gb \gep' n_{0}}
 \end{equation}           
       as soon as $\gb$ is large enough. 
       
     Indeed, there exists a vtj connected path of cycles $(\tilde{C}^{u}_{1}, \ldots,\tilde{C}^{u}_{m} )$ 
     of length $m$ (with  $ m \leq n_{0}$) joining $C(u)$ to $C(y)$. For any $u \in C$, applying the strong
      Markov property at the time of first entrance into $\tilde{C}^{u}_{1}$ and proceeding iteratively, we get:
       \begin{equation}
       \begin{aligned}
           \P_{\gb} \left[ C(y) \in (C^{u}_{1}, \ldots, C^{u}_{n_{0}})\right]  & \geq
           \P_{\gb} \left[ (C^{u}_{1}, \ldots, C^{u}_{m}) = (\tilde{C}^{u}_{1}, \ldots,\tilde{C}^{u}_{m} ) \right] \\
            & = \sum_{v \in \tilde{C}^{u}_{1} \cap \cB(C(u))} \P_{\gb} \left[ X^{u}_{\tau^{u}_{\cX \setminus C(u)}}=v, 
            (C^{u}_{2}, \ldots, C^{u}_{m}) = (\tilde{C}^{u}_{2}, \ldots,\tilde{C}^{u}_{m} ) \right] \\
             & \geq  e^{-\gb \gep'} \inf_{v \in \tilde{C}^{u}_{1}} \P_{\gb} \left[ 
            (C^{v}_{1}, \ldots, C^{v}_{m-1}) = (\tilde{C}^{u}_{2}, \ldots,\tilde{C}^{u}_{m} ) \right] \\
             & \geq \ldots \geq e^{-\gep' \gb n_{0}}
       \end{aligned}
       \end{equation} 
        where in the third inequality we used Corollary \ref{corex} and the definition of vtj--connectedness.  
      Since the last term does not depend on $u$, we get \eqref{unifdentroC}. 
     
%
       

       Making use recursively of the strong Markov property at times $\theta^{x}_{ke^{\gep \gb }/n_{0}}, k = 1,\ldots,n_{0}$,
       of the trivial bound $n_{0} \leq |\cX|$ and of \eqref{unifdentroC}, we get:  
    \begin{equation}
     \begin{aligned}
      \P_{\gb} \left[  j_{x,y} > e^{\gep \gb}  \right] & =  \P_{\gb} \left[ C(y) \notin (C^{x}_{1}, \ldots,C^{x}_{e^{\gep \gb}} ) \right]  \\
      & =  \sum_{\bar C \in \cM \setminus C(y)} \sum_{v \in C^{x}_{e^{\gb \gep} - n_{0}}}
      \P_{\gb} \left[ C(y) \notin (C^{x}_{1}, \ldots,C^{x}_{e^{\gb \gep} - n_{0}} ), X^{x}_{\theta_{e^{\gb \gep} - n_{0}}^{x}} = v,
      C^{x}_{e^{\gb \gep} - n_{0}} = \bar C \right] \\
     & \phantom{iiiiiiiiiiiiiiiiiiiiiiiiiiiiiiiiiiiiiiiii} \times  \P_{\gb} \left[ C(y) \notin (C^{v}_{1}, \ldots, C^{v}_{n_{0}})\right]  \\
      & \leq (1 - e^{-\gb \gep' n_{0}})  \P_{\gb} \left[ C(y) \notin (C^{x}_{1}, \ldots,C^{x}_{e^{\gb \gep} - n_{0}} ) \right] \\ 
       & \leq \ldots \\ 
      & \leq (1 - e^{- \gb \gep' |\cX|})^{e^{\gep \gb/|\cX|}}\\
       & \leq e^{- e^{\gb (\gep/|\cX| - \gep'|\cX|)} }
     \end{aligned}
    \end{equation}     
    and \eqref{ott2} then follows by choosing $\gep' \in (0,\gep/|\cX|^{2})$.      
This concludes the proof of Theorem \ref{crossgate}.
 
 \qed
\par\noindent

\medskip
\par\noindent
\textit{Proof of Theorem~\ref{crossgate}.\/}

 We first recall the following consequence of Proposition~\ref{decmaxcy}. 
  
  For any $i=1,\dots, n_0$:
   \begin{equation}\label{altsub}
    H(C_{i}) + \Gamma(C_{i}) = \Phi(x,y). 
   \end{equation} 
To get \eqref{altsub}, we first note 
that, since $y\in\cX\setminus C(x)$, by Propositions~\ref{OVirr} and 
\ref{fin00}, we have that 
$\Phi(x,y) \geq H(C(x)) + \Gamma(C(x))$. 
On the other hand, assume by contradiction that
$\Phi(x,y) >H(C(x))+\Gamma(C(x))$. Recalling 
\eqref{constsubcy} in Proposition~\ref{decmaxcy} and \eqref{usci1} 
in Proposition~\ref{usci}, it follows that there exists a path 
$\omega\in\Omega_{x,y}^{\rm{opt}}$ such that $\Phi_\omega<\Phi(x,y)$, which 
is absurd. 

Hence, we have proven that 
$\Phi(x,y)=H(C(x))+\Gamma(C(x))$.
By using again
\eqref{constsubcy} in Proposition~\ref{decmaxcy}, we then deduce \eqref{altsub}.

  Then we note that considering Proposition \ref{traje}, for Theorem \ref{crossgate} to hold, it is enough to show the inclusions 
  \begin{equation}\label{ottimali}
   \left( \Omega_{x,y}^{\rm{opt}} \cap E_{x,y} \right)
   \subset \cK_{x,y}  \subset \left\{ \tau^{x}_{W} < \tau^{x}_{y} \right\}.
  \end{equation}   
  
   Indeed, this implies in particular the trivial bound 
   \begin{displaymath}
   \P_{\gb}\left[\left\{ \tau^{x}_{W} < \tau^{x}_{y} \right\}\right] \geq \P_{\gb}\left[\cK_{x,y}\right],
   \end{displaymath} 
 and Proposition \ref{traje} provides the requested lower bound on this last quantity.

We remark that the 
inclusions of \eqref{ottimali} are strict in general.

The first inclusion follows immediately from the fact that 
an optimal path in $\Omega^\textrm{opt}_{x,y}$  exits
from an element of $\cM$ through its principal boundary. Also, it is clear that some paths 
in the set 
$\cK_{x,y}$ might
not be optimal, and hence that it might be strict in general. 

The second inclusion of \eqref{ottimali} is not straightforward 
and we stress that it 
relies crucially on Proposition \ref{usci}. Let us detail it. 
       
Consider first the case $\omega \in \cK_{x,y} \cap \Omega_{x,y}^{\rm{opt}}$.
Since $\omega \in \Omega_{x,y}^{\rm{opt}}$, by definition of a gate 
(see Section~\ref{s:sg}), it follows immediately that 
$\omega \cap W \neq \emptyset$.
          
Consider now  
an element $\omega \in \cK_{x,y} \setminus \Omega_{x,y}^{\rm{opt}}$, 
that is $\omega$ is an element 
of $\cK_{x,y}$
such that $\Phi(\omega) > \Phi(x,y)$.


 To show the second inclusion of \eqref{ottimali}, the strategy is the following: we consider the 
sequence of points $(u_1,\ldots,u_j)$ which are the successive points where $\omega$ 
intersects $\bigcup_{\bar{C} \in \cM} B(\bar{C})$. The sequence  $(u_1,\ldots,u_j)$ is nonempty from the construction
 of $\cK_{x,y}$ and from the fact that $C(x) \neq C(y)$. We are going to 
  construct stepwise a path $\widetilde{\omega} \in \cK_{x,y} \cap \Omega_{x,y}^{\rm{opt}}$ such that
 \begin{equation}\label{intertil}
  \widetilde{\omega} \bigcap \left(\bigcup_{\bar{C} \in \cM} B(\bar{C}) \right) = \{u_1,\ldots,u_j\}.
 \end{equation} 
 
  From the definition of a gate and from the fact that $\widetilde{\omega}$ is optimal, we deduce 
  that $\widetilde{\omega} \cap W \neq \emptyset$. From 
   this it follows that
    $\widetilde{\omega} \cap W = \omega \cap W \neq \emptyset$, which indeed implies
    the second inclusion of \eqref{ottimali}.

       To construct the path $\widetilde{\omega}$, we proceed in a recursive way;
       more precisely, we construct a sequence of paths $(\omega^{(k)})_{k \geq 0} \in \cK_{x,y}$ which 
        becomes stationary for $k$ large enough. We 
        initialize our recursion by setting $\omega^{(0)} := \omega$. Then, as long as the path $\omega^{(k)}$ is not 
        optimal, we proceed in the following way: consider 
        \begin{displaymath}
         i_{k} = \inf \left\{j \leq |\omega^{(k)}|, H\left(\omega^{(k)}_{j}\right) +
         \Delta\left(\omega^{(k)}_{j},\omega^{(k)}_{j+1}\right) > \Phi(x,y) \right\},
        \end{displaymath} 
         and 
         $C_{k}$ the element of $\cM$ containing $\omega^{(k)}_{i_{k}}$. Then we 
      distinguish two cases: $\omega^{(k)}_{i_{k}+1} \in \cB(C_{k})$ and $\omega^{(k)}_{i_{k}+1} \in C_{k}$.
       \begin{enumerate}
        \item[$\bullet$]  In the case where $\omega^{(k)}_{i_{k}+1} \in \cB(C_{k})$, we make use of \eqref{usci1} in Proposition
        \ref{usci} and of \eqref{altsub} to get that there exists 
         a path $\omega' \in \Omega_{\omega^{(k)}_{i_{k}},\omega^{(k)}_{i_{k}+1}}$ such that
         $\Phi(\omega') = \Gamma(C_{k}) + H(C_{k}) = \Phi(x,y)$ and
         for any $j \leq |\omega'|-1$,
          $\omega_{j}' \in C_{k}$. We define the concatenated path 
           \begin{equation}\label{concat}
            \omega^{(k+1)} := 
          \left( \left(\omega^{(k)}_{j}\right)_{j \leq i_{k}-1},\omega',\left(\omega^{(k)}_{j}\right)_{j \geq i_{k}+2}\right).
           \end{equation} 
           
            Note that 
          $\omega^{(k+1)} \in \cK_{x,y}$ and that $u \in \omega^{(k+1)}$. Then we continue the recursive construction. 
          
          \item[$\bullet$] In the case where $\omega^{(k)}_{i_{k}+1} \in C_{k}$, from \eqref{dent} in Proposition \ref{usci},
          there exists $ \omega' \in \Omega_{\omega^{(k)}_{i_{k}},\omega^{(k)}_{i_{k}+1}}$ 
        such that $\Phi(\omega') < \Phi(C_{k}, \cX \setminus C_{k})
       = \Phi(x,y)$ and such that $\omega'$ is entirely contained in $C_{k}$. Then we define the path $\omega^{(k+1)}$ as
       in \eqref{concat}, and we note that in this case
       also $\omega^{(k+1)} \in \cK_{x,y}$ and $u \in \omega^{(k+1)}$. 
       \end{enumerate}

        It is clear from the construction that the sequence of paths $(\omega^{(k)})_{k \geq 0}$ is stationary 
         after a number of steps at most $|\omega|$, and that the 
        final path $\widetilde{\omega}$ obtained at the end of the recursion  is an element
        of $\cK_{x,y} \cap \Omega_{x,y}^{\rm{opt}}$ satisfying \eqref{intertil}. Hence the second inclusion in \eqref{ottimali} follows, and thus 
        Theorem \ref{crossgate} is proved.

\qed
\par\noindent

Now we go to the proof of Proposition~\ref{Main}. We first note that, 
in the spirit of \cite{MNOS}, we need 
a downhill cycle path (see the 
definition in Section~\ref{tjsc}) connecting 
any given point $x\in\cX\setminus\cX_{a}$, for $a>0$, to 
$\cX_{a}$. 
We recall that the notion of downhill cycle path
given in \cite{MNOS} and \cite{OS1}, even if quite peculiar to the Metropolis 
dynamics setup, finds its natural extension to the general 
rare transition setup in \cite{OS2} and in \cite{CaCe} through the notion of 
"via typical jumps" connection. 
        
\medskip
\par\noindent
\textit{Proof of Proposition~\ref{Main}.\/}
Let $a>0$, 
we assume that $\cX_{a}$ is a proper subset of $\cX$, otherwise 
there is nothing to prove. 
We consider $x \in\cX\setminus\cX_{a}$ and note that, by 
Proposition~\ref{dhp}, there exists a vtj--connected 
cycle path $C_1,\ldots,C_l\subset\cX\setminus\cX_a$ such that
$x\in C_1$ and $\cB(C_l)\cap\cX_a\neq\emptyset$. 

Since none of the cycles $C_1,\ldots,C_l$
can contain points 
of $\cX_a$, for any $i=1,\ldots,l$ and any $z\in\cF(C_i)$ 
the stability level $V_z$ (recall definition \eqref{sl001}) 
of $z$ satisfies $V_{z} \leq a$, and hence from 
item~\ref{fin00i04} in Proposition~\ref{t:fin01}, we have
$\Gamma(C_i)\le a$ for any $i=1,\ldots,l$.

Then, from item~\ref{excyi01} in Proposition~\ref{excy},
for any cycle $C_i$ of the vtj--connected path,
for any $z\in C_i$, 
and for any $\gep > 0$, the function
\begin{displaymath}
\gb\in\R^{+}\mapsto\P_\beta
        \big[ \tau_{\partial C_i}^{z} > e^{\gb(a+ \gep)}\big]
\end{displaymath} 
is SES.

We consider $y\in\cB(C_l)\cap\cX_a$ and,
for each $2 \leq i \leq l$, we consider 
$y_{i} \in \cB(C_{i-1}) \cap C_{i}$. We 
define $y_{1} = x$
and
$y_{l+1} = y$, and we consider the set of paths
\begin{equation}\label{defE}
\cE :=
\cE\big((C_{1},x),(C_{2},y_{2}),\ldots,(C_{l},y_{l}),(\cX_{a},y) \big)
\end{equation}
consisting of the paths constructed by the concatenation 
of any $l$--uple of paths 
$\go^{1},\go^{2},\ldots,\go^{l}$ satisfying the following conditions:
\begin{enumerate}
\item  \label{c1}
for any $i=1.\ldots,l$ the length of the path 
$\go^{i}$ satisfies $|\go_{i}| \leq e^{\gb(a+\gep/4)}$;
\item 
for any $i=1,\ldots,l$ the path $\go^i$ joins $y_{i}$ to $y_{i+1}$, 
that is, $\go^{i} \in \Omega_{y_{i},y_{i+1}}$ (recall the notation 
introduced in Section~\ref{s:defpre});
\item 
$\go^{i}_{j} \in C_{i}$
for any $i=1,\ldots,l$ and for any 
$j=1,\ldots,|\go^{i}| -1$.
\end{enumerate}
 
The existence of such a family of paths is ensured by 
Propositions~\ref{OVirr}, \ref{fin00}, \ref{usc0}, and \ref{usci}.  
We stress that condition \ref{c1} restricts the set $\cE$ to 
 paths which spend a time less than 
$e^{\gb(a + \gep/4)}$ in any cycle $C_{i}, i \leq l$.
        
For shortness, in the sequel, we shall use the  notation 
$\cE$ for 
the set of trajectories defined in \eqref{defE}. 
 
        
Note that the length of any $\go \in \cE$ 
satisfies the upper bound 
$|\go| \leq |\cX| e^{\gb(a+\gep/4)}$. Moreover, since the state space $\cX$ 
is finite, we can assume that $\gb$ is large enough so that 
\begin{displaymath}
|\go| \leq |\cX|\, e^{\gb(a+\gep/4)}
\le e^{\gb(a+\gep/2)}
\;\;\textrm{ for any }
\go\in\cE.
\end{displaymath}

Now, we write 
\begin{displaymath}
\P_\gb\big[ \tau_{\cX_{a}}^{x} \leq e^{\gb(a + \gep/2)}\big]
\geq 
\P_\gb\big[\tau_{\cX_{a}}^{x} 
  \leq e^{\gb(a + \gep/2)},(X_{k})_{k \leq \tau^x_{\cX_{a}}} 
\in \cE \big] 
= 
\P_\gb\big[(X_{k})_{k \leq \tau^x_{\cX_{a}}} \in \cE \big],
\end{displaymath} 
where in the last step we have used the bound above
on the length of the trajectories in $\cE$. 
Then we use Markov's property to get that
\begin{displaymath}
\P_\gb\big[ \tau_{\cX_{a}}^{x} \leq e^{\gb(a + \gep/2)}\big]
\geq 
\P_\gb\big[
\left(X_{k}\right)_{k \leq \tau_{\cX_{a}}^{x}} 
\in \cE \big]
= 
\prod_{i=1}^{l} 
\P_\gb\big[ 
   \tau^{y_{i}}_{\cX\setminus C_{i}} 
   \leq e^{\gb(a + \gep/4)}, X^{y_{i}}_{\tau^{y_i}_{\cX\setminus C_{i}}}
= y_{i+1} \big].
\end{displaymath} 
 
Combining this inequality and \eqref{triv} implies that, for any $\gep'>0$,
\begin{displaymath}
\P_\gb\big[\tau_{\cX_{a}}^{x} \leq e^{\gb(a + \gep/2)}\big]
\geq 
e^{- \gb \gep' l} 
\geq e^{- \gb \gep' |\cX|}
\end{displaymath}
as soon as $\gb$ is large enough.
        
Since  the last term in the right hand side of the bound above 
does not depend on $x \in \cX_{a}$, we get that 
\begin{displaymath}
\inf_{x \in \cX_{a}} 
\P_\gb\big[ \tau_{\cX_{a}}^{x} \leq e^{\gb(a + \gep/2)}\big]
\geq 
e^{- \gb \gep' |\cX|}.
\end{displaymath}
   
Now we iterate this inequality by making use of the Markov's property at 
the times $ke^{\gb(a + \gep/2)}$, $k=1,\ldots, e^{\gb \gep /2}$
to get that
\begin{displaymath}
\P_\gb\big[\tau_{\cX_{a}}^{x} > e^{\gb(a + \gep)}\big] 
\leq  
\Big(\sup_{x'\in\cX_{a}}
 \P_\gb\big[\tau_{\cX_{a}}^{x'} > 
    e^{\gb(a + \gep/2)}\big]\Big)^{e^{\gb \gep /2}}
\!\!
\!\!
\leq  
\Big( 1 - e^{-\gb \gep' |\cX|} \Big)^{e^{\gb \gep /2}} 
\!\!
\!\!
\leq 
e^{ - e^{\gb (\gep /2-\gep' |\cX|) }}
\end{displaymath} 
for any $x \in \cX_{a}$.
   
   Finally, picking up $\gep' > 0$ small enough, we get that the function $\gb \mapsto e^{ - e^{ \gb (\gep /2-\gep' |\cX|) }}$ is SES, and thus 
Proposition~\ref{Main} is proved.
\qed

\medskip
\par\noindent
\textit{Proof of Proposition~\ref{t:staz03}.\/}
Set
$T:=\exp \left( \beta(h-\varepsilon) \right)$; 
writing $x_0=x$ and making use of the Markov property, we immediately get: 
\begin{displaymath}
\P_\gb(\cE^{x,h}(\varepsilon))
\le
   \sum_{n=1}^{\lfloor T\rfloor} \P_\gb(\cE_n^{x,h})
   =
   \sum_{n=1}^{\lfloor T\rfloor}
   \sum_{\newatop{x_1,\dots,x_n\in \cX:}
                 {H(x_{n-1})+\Delta(x_{n-1},x_n)\ge H(x)+h}} 
   p_\beta(x,x_1)\cdots p_\beta(x_{n-1},x_n).
\end{displaymath}
 
We multiply and divide by $\mu_\beta(x)$ on the right hand side (recall 
that for $\gb$ large enough, the Markov chain 
is irreducible, and hence $\mu_{\gb}$ is strictly positive over $\cX$). 
Also, we estimate the first two terms with the sum over the first state, and 
we deduce
\begin{displaymath}
\P_\gb(\cE^{x,h}(\varepsilon))
\le
   \frac{1}{\mu_\beta(x)}
   \sum_{n=1}^{\lfloor T\rfloor}
   \!\!
   \!\!
   \sum_{\newatop{x_1,\dots,x_n\in \cX:}
                 {\newatop{H(x_{n-1})+\Delta(x_{n-1},x_n)}{\ge H(x)+h}}} 
   \!\!
   \!\!
   \!\!
   \Big[
   \sum_{x_0\in \cX}
   \mu_\beta(x_0)
   p_\beta(x_0,x_1)
   \Big]
   p_\beta(x_1,x_2)\cdots p_\beta(x_{n-1},x_n)
\end{displaymath}
 
Now, 
making use of the stationarity of $\mu_{\gb}$, we get  
\begin{displaymath}
\P_\gb(\cE^{x,h}(\varepsilon))
\le
   \frac{1}{\mu_\beta(x)}
   \sum_{n=1}^{\lfloor T\rfloor}
   \sum_{\newatop{x_{n-1},x_n\in \cX:}
                 {H(x_{n-1})+\Delta(x_{n-1},x_n)\ge H(x)+h}} 
   \mu_\beta(x_{n-1})p_\beta(x_{n-1},x_n),
\end{displaymath}
and hence
\begin{displaymath}
\P_\gb(\cE^{x,h}(\varepsilon))
\le
   \frac{\lfloor T\rfloor\,|\cX|^2}{\mu_\beta(x)}
   \sup_{\newatop{w,z\in \cX:}
                 {H(w)+\Delta(w,z)\ge H(x)+h}} 
   \mu_\beta(w)p_\beta(w,z)
\end{displaymath}
 
Recalling the convergences 
\eqref{FW} and \eqref{virth}, we get that for any $\gep' > 0$, as soon as $\gb$ is large enough:
\begin{align*}
\P_\gb(\cE^{x,h}(\varepsilon))
& \le
\lfloor T\rfloor|\cX|^2 e^{\beta (H(x)+ \gep')}
   \sup_{\newatop{w,z\in \cX:}
                 {H(w)+\Delta(w,z)\ge H(x)+h}} 
    e^{-\beta[H(w)+\Delta(w,z)- \gep']} \\ 
    & \le
 \lfloor T\rfloor |\cX|^2
  e^{\beta (H(x)+ 2 \gep')} e^{-\beta(H(x)+h)} \\
   & \le
  |\cX|^2
  e^{\beta ( 2 \gep' - \gep)},
\end{align*}
where in the last step we used the definition of $T$. Now, choosing $\gep' \in (0,\gep/4)$ concludes the proof of Proposition~\ref{t:staz03}.
\qed

\medskip
\par\noindent
\textit{Proof of Proposition~\ref{t:quasimet}.\/}
Consider $ C \in \cC(\cX)$, $z \in \cF(C)$ and $\gep > 0$. By the finiteness of $\cF(C)$, it is enough to prove 
\begin{equation}\label{qumetw}
\P_{\gb}\left[\Phi((X^{z}_{t})_{0\leq t 
      \leq \tau^{z}_{\cX \setminus C}}) >H(C) + \Gamma(C) + \gep \right] 
\leq e^{-\gb \gep/4} 
\end{equation} 
for \eqref{quasimet} to hold. 
We consider the following decomposition:
\begin{equation}
\label{dec000}
\begin{array}{l}
\P_{\gb}\left[\Phi((X^{z}_{t})_{0\leq t \leq \tau^{z}_{\cX \setminus C}}) 
     > H(C) + \Gamma(C) + \gep \right]
\vphantom{\bigg\{_\}}\\
\phantom{mmmmmmm}=  
    \P_{\gb}\left[\Phi((X^{z}_{t})_{0\leq t \leq \tau^{z}_{\cX \setminus C}}) 
    > H(C) + \Gamma(C) + \gep, \tau^{z}_{\cX \setminus C} \leq e^{\gb(\Gamma(C) + \gep/2)}  \right] 
\vphantom{\bigg\{_\}}\\
\phantom{mmmmmmm=}  
     + \P_{\gb}\left[\Phi((X^{z}_{t})_{0\leq t \leq \tau^{z}_{\cX \setminus C}}) 
    > H(C) + \Gamma(C) + \gep,
     \tau^{z}_{\cX \setminus C} > e^{\gb(\Gamma(C) + \gep/2)}  \right]. 
    \end{array}
   \end{equation} 
   
For the second term in the right hand side above, 
we deduce from item~\ref{excyi01} in Proposition~\ref{excy} that
\begin{equation}\label{fin01}
\begin{array}{l}
\P_{\gb}\left[\Phi((X^{z}_{t})_{0\leq t \leq \tau^{z}_{\cX \setminus C}}) > H(C) + \Gamma(C) + \gep,
     \tau^{z}_{\cX \setminus C} > e^{\gb(\Gamma(C) + \gep/2)}  \right] \\
\phantom{iiiiiiiiiiiiiiiiiiiiiiiiiiiiiiiiiiiiiiiiiiiiiiii}
\leq
     \P_{\gb}\left[
     \tau^{z}_{\cX \setminus C} > e^{\gb(\Gamma(C) + \gep/2)}  \right] 
       \leq e^{-\gb \gep/4}
     \end{array}
    \end{equation} 
  as $\gb \to \infty$. 

As for the first term, we first have the inequality:
   \begin{equation}\label{fin02}
   \begin{aligned}
     &\P_{\gb}\left[\Phi((X^{z}_{t})_{0\leq t \leq \tau^{z}_{\cX \setminus C}}) > H(C) + \Gamma(C) + \gep,
     \tau^{z}_{\cX \setminus C} \leq e^{\gb(\Gamma(C) + \gep/2)}  \right] \\  & \phantom{iiiiiiiiiiiiiiiiiiiiiiiiiiiiiii}  \leq
     \P_{\gb}\left[\Phi((X^{z}_{t})_{0\leq t \leq e^{\gb(\Gamma(C)+\gep/2)}})  > H(C) + \Gamma(C) + \gep \right]. 
   \end{aligned}
   \end{equation}  
   Using the fact that $z \in \cF(C)$ (that is $H(z) = H(C)$), we have the equality of events:
   \begin{displaymath}
    \left\{\Phi((X^{z}_{t})_{t \leq e^{\gb(\Gamma(C)+\gep/2)}}) > H(C) +  \Gamma(C) + \gep \right\} = \cE^{z,\Gamma(C) + \gep}(\gep/2),
   \end{displaymath} 
where we have recalled \eqref{staz04} and \eqref{staz05}.
We deduce from Proposition~\ref{t:staz03} that the term in the right hand side of \eqref{fin02} is less than $e^{-\gb \gep/4}$ as $\gb \to \infty$. 
Combining this
inequality 
with 
\eqref{dec000} and 
\eqref{fin01}, we deduce \eqref{qumetw}, 
and, hence, Proposition~\ref{t:quasimet}. 
\qed
           
\appendix

\section{Computing differences of virtual energy}
\label{s:gener}
In this appendix, we describe an abstract framework for which 
the virtual energy has \textit{a priori} no explicit expression, but where
we can construct it stepwise starting from a reference point acting as a point of null potential. 
 
We consider a Freidlin Wentzell dynamics satisfying Definition~\ref{setupFW} 
and such that for every $x,y\in \cX$
   \begin{equation}\label{symdy}
    \Delta(x,y) < \infty \hspace{6 pt} \text{if and only if} \hspace{6 pt}    \Delta(y,x) < \infty.        
   \end{equation}
            
      Moreover, we  assume that the dynamics  
   satisfies the additional condition (where we recall that $\mu_{\gb}$ is 
the invariant measure). 
   
 \begin{condition}
For any $\beta>0$, 
there exists a function
$\rho:\R_+\to\R_+$ such that
$\rho(\beta) \to 0$ as $\beta\to\infty$ and 
\begin{equation}
\label{qr}
\Big|
-\log\mu_\beta(x)+\beta\Delta(x,y)
-\left[-\log\mu_\beta(y)+\beta\Delta(y,x)\right]
\Big|
\le\beta\rho(\beta)
\end{equation}
for any $x,y\in\cX$.
\end{condition}

 Of course, the convergence \eqref{qr} is nothing else than requesting the existence of a potential, which is equal to the virtual
  energy up to a constant (see \eqref{weakrev} and Proposition \ref{t:ve020}).

  Now we fix an arbitrary state $\bar x \in \cX$ and we define the Hamiltonian--like quantity
\begin{equation}
\label{defg}
G_\beta(x):=-\log[\mu_\beta(x)/\mu_\beta(\bar x)].
\end{equation}
  
   For any $x \in \cX, x \neq \bar x$, by irreducibility, there exists 
a path $\go\in\Omega_{\bar x, x}$ such that $|\go| \leq \cX$.
    Given such a path, we define the quantity 
    \begin{equation}\label{defW}
     W_{\go}(x):= \sum_{i=2}^{|\go|} \left[\Delta(\go_{i-1},\go_i)-\Delta(\go_i,\go_{i-1}) \right]
    \end{equation} 
   and we set $W_{\go}(\bar x):= 0$.
  
  \begin{proposition}
 Given $x\in\cX$ and $x\neq\bar x$, 
the quantity $W_{\go}(x)$ defined by \eqref{defW}  
does not depend on the particular 
choice of the path $\omega \in \Omega_{\bar x,x}$, and hence it 
defines a function $W: \cX \to \R$. The function $W(\cdot) - \min_{\cX} W$ coincides with the virtual energy $H$.
  \end{proposition}
 
  In general, the virtual energy might have an expression too involved for practical purposes. Equation \eqref{defW} 
   provides a constructive way to compute explicitly $H$ step by step just from the knowledge of the rates of the dynamics.

\smallskip
\par\noindent
\textit{Proof.\/} 
For any $x,y\in\cX$ and $x\neq y$, 
we consider $\go,\go' \in \Omega_{x,y}$ and show that 
    \begin{equation}\label{ompri}
     \sum_{i=2}^{|\go|} \left[\Delta(\go_{i-1},\go_i)-\Delta(\go_i,\go_{i-1}) \right] = \sum_{i=2}^{|\go'|}
     \left[\Delta(\go'_{i-1},\go'_i)-\Delta(\go'_i,\go'_{i-1}) \right].
    \end{equation}  
    
Indeed,  
using telescoping sums in the right hand side above, 
we can assume that all the $\go'_{i}$'s are distinct (and in 
particular $|\go'| \leq |\cX|$). 
   
    By \eqref{qr}, we get the inequality
    \begin{displaymath}
     \begin{aligned}
      \big|
 & G_{\gb}(x)-\gb \sum_{i=2}^{|\go|}
 [\Delta(\go_{i-1},\go_i)-\Delta(\go_i,\go_{i-1})]
 \big|
\\ & \phantom{.} =
 \big|
 \sum_{i=2}^{|\go|}[G_{\gb}(\go_i)-G_{\gb}(\go_{i-1})]
 -\gb
 \sum_{i=2}^{|\go|}[\Delta(\go_{i-1},\go_i)-\Delta(\go_i,\go_{i-1})]
 \big|
 \leq|\cX|\beta\rho(\gb).      
     \end{aligned}
    \end{displaymath} 
    
     By triangular inequality, we then deduce that
  \begin{displaymath}
   \begin{aligned}
   & \big| \gb \sum_{i=2}^{|\go|} [\Delta(\go_{i-1},\go_i)-\Delta(\go_i,\go_{i-1})] 
    - \gb \sum_{i=2}^{|\go'|} [\Delta(\go'_{i-1},\go'_i)-\Delta(\go'_i,\go'_{i-1})] \big|
   \\ & \phantom{iii} = \big| G_{\gb}(x) - \gb \sum_{i=2}^{|\go|} [\Delta(\go_{i-1},\go_i)-\Delta(\go_i,\go_{i-1})]   \big|  
\\ &
\phantom{iii=}
+ 
     \big| G_{\gb}(x) - \gb \sum_{i=2}^{|\go'|} [\Delta(\go'_{i-1},\go'_i)-\Delta(\go'_i,\go'_{i-1})] \big| \\
      & \phantom{iii} \leq  2|\cX|\beta\rho(\gb).
   \end{aligned}
  \end{displaymath}     
Now we divide both sides by $\gb$ and we let $\gb \to \infty$ 
to deduce \eqref{ompri}.

\qed
  
\section{Explicit expression of the virtual energy}
\label{s:graph}
As noted in Section~\ref{virth}, the virtual energy $H(x)$, 
for $x\in\cX$, has 
an explicit expression in terms of a specific graph construction. 
The same holds for the functions $\Gamma_D(x)$ and $\Delta_D(x,y)$,
with $D\subset\cX$, $x\in D$, and $y\in\cX\setminus D$, introduced 
in Proposition~\ref{exgen}.
These explicit expressions 
were not necessary for our purposes, but for the sake of completeness, we choose to summarize these formulas in 
this appendix. 

We use the notations of \cite{Ca}, but since we do not want to 
develop the full theory here, we try to keep it as minimal 
as possible. 

\begin{definition}
\label{graph000}
Given $A\subset\cX$ nonempty, let $G(A)$ be the set of 
\emph{oriented graphs} $g\in\cX\times\cX$ verifying the following properties: 
\begin{itemize}
\item[--]
for any $x\in\cX\setminus A$, there exists a unique $y\in\cX$ 
such that $(x,y)\in g$ (namely for any point in $\cX\setminus A$, there 
exists a unique arrow of the graph $g$ exiting from such a point);
\item[--]
for any edge $(x,y)\in g$, $x\in\cX\setminus A$ (no arrow of the graph 
$g$ exits from $A$);
\item[--]
for any $x\in\cX$, $n\in\N$, 
$(x,x_1),(x_1,x_2),\dots,(x_{n-1},x_n)\in g$ 
one has that $x\neq x_i$ for $i=1,\dots,n$
(the graph $g$ is without loops).
\end{itemize}
\end{definition}

Since $\cX$ is finite, 
from this definition it follows 
that for $x\in\cX\setminus A$, there exists a sequence of 
arrows connecting $x$ to $A$.
We borrow (and adapt to our notation) 
a beautiful description of the set $G(A)$ 
from \cite[below Definition~3.1]{OS2}: $G(A)$ is a 
forest of trees with roots in $A$ and with branches given 
by arrows directed towards the root. 

\begin{definition}
\label{graph010}
Given $A\subset\cX$ nonempty, $x\in\cX\setminus A$, 
and $y\in A$, let $G_{x,y}(A)$ 
be the collection of graphs $g\in G(A)$ such that 
there exist $n\in\N$ and 
$x_1,\dots,x_n\in\cX$ such that 
$(x,x_1),(x_1,x_2)\dots,(x_n,y)\in g$. 
\end{definition}

In words, $G_{x,y}(A)$ is the set of graphs in $G(A)$
 connecting the point $x$ to the point $y$. 

   For any $x\in\cX$, the virtual energy $H(x)$ is given by (see \cite[Proposition~4.1]{Ca})
\begin{displaymath}
H(x)=
\min_{g\in G(\{x\})}\sum_{(w,z)\in g}\Delta(w,z)
-
\min_{x'\in\cX}
\min_{g\in G(\{x'\})}\sum_{(w,z)\in g}\Delta(w,z).
\end{displaymath}
 
Moreover (see \cite[Proposition~4.2]{Ca}), 
for any 
$D\subset\cX$ nonempty, $x\in D$, and $y\in\cX\setminus D$, one has the following equality: 
\begin{displaymath}
\Gamma_D(x)
=
\min_{g\in G(\cX\setminus D)}\sum_{(w,z)\in g}\Delta(w,z)
-
\min_{x'\in\cX\setminus D}
\min_{g\in G_{x,x'}((\cX\setminus D)\cup\{x'\})}\sum_{(w,z)\in g}\Delta(w,z),
\end{displaymath}
and similarly  
\begin{displaymath}
\Delta_D(x,y)
=
\min_{g\in G_{x,y}(\cX\setminus D)}\sum_{(w,z)\in g}\Delta(w,z)
-
\min_{g\in G(\cX\setminus D)}\sum_{(w,z)\in g}\Delta(w,z).
\end{displaymath}

\bigskip

{\bf Acknowledgments.\/}  JS acknowledges  
the support of NWO STAR grant: "Metastable and cut--off 
behavior of stochastic processes" and 
thanks University of Roma Tre for the kind hospitality. The authors thank 
A.\ Zocca for useful comments and discussions. 
ENMC thanks the Institute of Complex Molecular Systems (TU/e, Eindhoven)
for the kind hospitality.

\bibliographystyle{plain}

\bibliography{IrrSetbibli}

\begin{thebibliography}{10}

\bibitem{AB1}
D.~J. Aldous and M.~Brown.
\newblock Inequalities for rare events in time-reversible {M}arkov chains. {I}.
\newblock In {\em Stochastic inequalities ({S}eattle, {WA}, 1991)}, volume~22
  of {\em IMS Lecture Notes Monogr. Ser.}, pages 1--16. Inst. Math. Statist.,
  Hayward, CA, 1992.

\bibitem{AB2}
D.~J. Aldous and M.~Brown.
\newblock Inequalities for rare events in time-reversible {M}arkov chains.
  {II}.
\newblock {\em Stochastic Process. Appl.}, 44(1):15--25, 1993.

\bibitem{BL1}
J.~Beltr{\'a}n and C.~Landim.
\newblock Metastability of reversible finite state {M}arkov processes.
\newblock {\em Stochastic Process. Appl.}, 121(8):1633--1677, 2011.

\bibitem{BL2}
J.~Beltr{\'a}n and C.~Landim.
\newblock Tunneling and metastability of continuous time {M}arkov chains {II},
  the nonreversible case.
\newblock {\em J. Stat. Phys.}, 149(4):598--618, 2012.

\bibitem{BLM}
O.~Benois, C.~Landim, and M.~Mourragui.
\newblock Hitting times of rare events in {M}arkov chains.
\newblock {\em J. Stat. Phys.}, 153(6):967--990, 2013.

\bibitem{BEGK04}
Anton Bovier, Michael Eckhoff, V{\'e}ronique Gayrard, and Markus Klein.
\newblock Metastability in reversible diffusion processes. {I}. {S}harp
  asymptotics for capacities and exit times.
\newblock {\em J. Eur. Math. Soc. (JEMS)}, 6(4):399--424, 2004.

\bibitem{CGOV}
Marzio Cassandro, Antonio Galves, Enzo Olivieri, and Maria~Eul{\'a}lia Vares.
\newblock Metastable behavior of stochastic dynamics: a pathwise approach.
\newblock {\em J. Statist. Phys.}, 35(5-6):603--634, 1984.

\bibitem{Ca}
O.~Catoni.
\newblock Simulated annealing algorithms and {M}arkov chains with rare
  transitions.
\newblock In {\em S\'eminaire de {P}robabilit\'es, {XXXIII}}, volume 1709 of
  {\em Lecture Notes in Math.}, pages 69--119. Springer, Berlin, 1999.

\bibitem{CaCe}
Olivier Catoni and Rapha{\"e}l Cerf.
\newblock The exit path of a {M}arkov chain with rare transitions.
\newblock {\em ESAIM Probab. Statist.}, 1:95--144 (electronic), 1995/97.

\bibitem{CN03}
E.~N.~M. Cirillo and F.~R. Nardi.
\newblock Metastability for a stochastic dynamics with a parallel heat bath
  updating rule.
\newblock {\em J. Statist. Phys.}, 110(1-2):183--217, 2003.

\bibitem{CNSo}
E.~N.~M. Cirillo, F.~R. Nardi, and J.~Sohier.
\newblock A comparison between different cycle decompositions for metropolis
  dynamics.
\newblock {\em Markov Process. Related Fields}, To appear.

\bibitem{CNS08p}
E.~N.~M. Cirillo, F.~R. Nardi, and C.~Spitoni.
\newblock Competitive nucleation in reversible probabilistic cellular automata.
\newblock {\em Phys. Rev. E (3)}, 78(4):040601, 4, 2008.

\bibitem{CNS08j}
E.~N.~M. Cirillo, F.~R. Nardi, and C.~Spitoni.
\newblock Metastability for reversible probabilistic cellular automata with
  self-interaction.
\newblock {\em J. Stat. Phys.}, 132(3):431--471, 2008.

\bibitem{CLRS}
Emilio N.~M. Cirillo, Pierre-Yves Louis, Wioletta~M. Ruszel, and Cristian
  Spitoni.
\newblock {Effect of self-interaction on the phase diagram of a Gibbs-like
  measure derived by a reversible Probabilistic Cellular Automata}.
\newblock {\em {Chaos Solitons \& Fractals}}, {64}({SI}):{36--47}, {JUL}
  {2014}.

\bibitem{CNP}
E.N.M. Cirillo, F.R. Nardi, and A.D. Polosa.
\newblock {Magnetic order in the Ising model with parallel dynamics}.
\newblock {\em {Physical Review E}}, {64}({5, 2}), {2001}.

\bibitem{FMNS}
R.~Fernandez, F.~Manzo, F.~R. Nardi, and E.~Scoppola.
\newblock Asymptotically exponential hitting times and metastability: a
  pathwise approach without reversibility.
\newblock {\em http://arxiv.org/abs/1406.2637}.

\bibitem{FMNSS}
R.~Fernandez, F.~Manzo, F.~R. Nardi, E.~Scoppola, and J.~Sohier.
\newblock Conditioned, quasi-stationary, restricted measures and metastability.
\newblock {\em http://arxiv.org/abs/1234.5678}.

\bibitem{GJH}
G.~Grinstein, C.~Jayaprakash, and Yu~He.
\newblock Statistical mechanics of probabilistic cellular automata.
\newblock {\em Phys. Rev. Lett.}, 55(23):2527--2530, 1985.

\bibitem{MNOS}
F.~Manzo, F.~R. Nardi, E.~Olivieri, and E.~Scoppola.
\newblock On the essential features of metastability: tunnelling time and
  critical configurations.
\newblock {\em J. Statist. Phys.}, 115(1-2):591--642, 2004.

\bibitem{MRRTT}
N.~Metropolis, A.W. Rosenbluth, M.N. Rosenbluth, A.H. Teller, and E.~Teller.
\newblock {Equation of state calculations by fast computing machines.}
\newblock {\em {Journal of Chemical Physics}}, {21}({6}):{1087--1092}, {1953}.

\bibitem{OI}
R.~I. Oliveira.
\newblock Mean field conditions for coalescing random walks.
\newblock {\em Ann. Probab.}, 41(5):3420--3461, 2013.

\bibitem{OS1}
E.~Olivieri and E.~Scoppola.
\newblock Markov chains with exponentially small transition probabilities:
  first exit problem from a general domain. {I}. {T}he reversible case.
\newblock {\em J. Statist. Phys.}, 79(3-4):613--647, 1995.

\bibitem{OS2}
E.~Olivieri and E.~Scoppola.
\newblock Markov chains with exponentially small transition probabilities:
  first exit problem from a general domain. {II}. {T}he general case.
\newblock {\em J. Statist. Phys.}, 84(5-6):987--1041, 1996.

\bibitem{OV}
E.~Olivieri and M.~E. Vares.
\newblock {\em Large deviations and metastability}, volume 100 of {\em
  Encyclopedia of Mathematics and its Applications}.
\newblock Cambridge University Press, Cambridge, 2005.

\bibitem{Tr}
A.~Trouv{\'e}.
\newblock Cycle decompositions and simulated annealing.
\newblock {\em SIAM J. Control Optim.}, 34(3):966--986, 1996.

\end{thebibliography}
 
\end {document}